\documentclass[12pt,reqno,twoside]{amsart}
\usepackage{amsfonts}
\usepackage{amsmath,amssymb,amsthm,amsthm}
\usepackage{mathtools}
\usepackage{dsfont}
\usepackage{etoolbox}
\usepackage{color}
\usepackage{algorithm}
\usepackage{algorithmic}
\usepackage{hyperref}
\usepackage{url}

\usepackage{enumitem}
\setlist[enumerate]{leftmargin=.5in}
\setlist[itemize]{leftmargin=.5in}

\floatname{algorithm}{Algorithm}

\usepackage[dvips,bottom=1.4in,right=1in,top=1in, left=1in]{geometry}
\usepackage[capitalize,nameinlink]{cleveref}[0.19]
\crefname{section}{section}{sections}
\crefname{subsection}{subsection}{subsections}
\Crefname{section}{Section}{Sections}
\Crefname{subsection}{Subsection}{Subsections}
\crefname{hypothesis}{Hypothesis}{Hypotheses}
\Crefname{figure}{Figure}{Figures}
\crefformat{equation}{\textup{#2(#1)#3}}
\crefrangeformat{equation}{\textup{#3(#1)#4--#5(#2)#6}}
\crefmultiformat{equation}{\textup{#2(#1)#3}}{ and \textup{#2(#1)#3}}
{, \textup{#2(#1)#3}}{, and \textup{#2(#1)#3}}
\crefrangemultiformat{equation}{\textup{#3(#1)#4--#5(#2)#6}}%
{ and \textup{#3(#1)#4--#5(#2)#6}}{, \textup{#3(#1)#4--#5(#2)#6}}{, and \textup{#3(#1)#4--#5(#2)#6}}

\Crefformat{equation}{#2Equation~\textup{(#1)}#3}
\Crefrangeformat{equation}{Equations~\textup{#3(#1)#4--#5(#2)#6}}
\Crefmultiformat{equation}{Equations~\textup{#2(#1)#3}}{ and \textup{#2(#1)#3}}
{, \textup{#2(#1)#3}}{, and \textup{#2(#1)#3}}
\Crefrangemultiformat{equation}{Equations~\textup{#3(#1)#4--#5(#2)#6}}%
{ and \textup{#3(#1)#4--#5(#2)#6}}{, \textup{#3(#1)#4--#5(#2)#6}}{, and \textup{#3(#1)#4--#5(#2)#6}}

\crefdefaultlabelformat{#2\textup{#1}#3}

\usepackage[textsize=small]{todonotes}
\def\RR{{\mathbb{R}}}
\def\NN{{\mathbb{N}}}
\def\R{{\mathcal{R}}}

\newtheorem{theorem}{Theorem}[section]
\newtheorem{lemma}[theorem]{Lemma}
\newtheorem{remark}[theorem]{Remark}
\newtheorem{example}[theorem]{Example}
\newtheorem{definition}[theorem]{Definition}

\title{Fractional Deep Neural Network via Constrained Optimization}

\author{Harbir Antil$^1$, Ratna Khatri$^1$, Rainald L\"ohner$^2$, and Deepanshu Verma$^1$} 
\address{$^1$Department of Mathematical Sciences and the Center for Mathematics and Artificial Intelligence (CMAI), George Mason University, Fairfax, VA 22030, USA.} 
\address{$^2$Department of Computational and Data Science, George Mason University, Fairfax, Virginia, USA.}
\email{hantil@gmu.edu, rkhatri3@gmu.edu, rlohner@gmu.edu, dverma2@gmu.edu}

\thanks{This work is partially supported by Department of Navy, Naval Postgraduate School - N00244-20-1-0005, Air Force Office of Scientific Research (AFOSR) under Award NO: FA9550-19-1-0036, and National Science Foundation grants DMS-1818772 and DMS-1913004. The second author is also partially supported by a Provost award at George Mason University under the Industrial Immersion Program.}



\begin{document}

\begin{abstract}
This paper introduces a novel algorithmic framework for a deep neural network (DNN),  which in a mathematically rigorous manner, allows us to incorporate history (or memory) into the network -- it ensures all layers are connected to one another. This DNN, called Fractional-DNN, can be viewed as a time-discretization of a fractional in time nonlinear ordinary differential equation (ODE). The learning problem then is a minimization problem subject to that fractional ODE as constraints. We emphasize that an analogy between the existing DNN and ODEs, with standard time derivative, is well-known by now. The focus of our work is the Fractional-DNN. Using the Lagrangian approach, we provide a derivation of the backward propagation and the design equations. We test our network on several datasets for classification problems. Fractional-DNN offers various  advantages over the existing DNN. The key benefits are a significant improvement to the vanishing gradient issue due to the memory effect, and better handling of nonsmooth data due to the network's ability to approximate non-smooth functions.
\end{abstract}
\maketitle

\section{Introduction}
Deep learning has emerged as a potent area of research and has enabled a remarkable progress in recent years spanning domains like 
imaging science \cite{he2016deep, ADK_19_ml_tomo, wu2018deep_imaging, CNN_tomo_2017}, biomedical applications \cite{lee2018deepres_bio, chen2018voxresnet_bio, Pock2017_MRI}, satellite imagery, remote sensing \cite{Tai_2017_CVPR, zhang2018missing, bischke2017detection}, etc.
However, the mathematical foundations of many machine learning architectures are largely lacking \cite{Bengio2010, Qiu2016, Wigderson2019, ruthotto2018deep, Weinan_2019}. The current trend of success is largely due to the empirical evidence. Due to the lack of mathematical foundation, it becomes challenging to understand the detailed workings of networks \cite{goldt2019modelling, mallat2013deep}.

The overarching goal of machine learning algorithms is to learn a function using some known data. Deep Neural Networks (DNN), like Residual Neural Networks (RNN), are a  popular family of deep learning architectures which have turned out to be groundbreaking in imaging science. An introductory example of RNN is the ResNet \cite{he2016deep} which has been successful for classification problems in imaging science. Compared to the classical DNNs, the innovation of the RNN architecture comes from a simple addition of an identity map between each layer of the network. This ensures a continued flow of information from one layer to another. Despite their success, DNNs are prone to various challenges such as vanishing gradients \cite{Bengio1994,  Bengio2010, veit2016residual}, difficulty in approximating non-smooth functions, long training time \cite{Haber2017_RNN}, etc.
  
We remark that recently in \cite{huang2017denseNet} the authors have introduced a DenseNet, which is a new approach to prevent the gradient ``wash out" by considering dense blocks, in which each layer takes into account all the previous layers (or the memory). They proceed by concatenating the outputs of each dense block which is then fed as an input to the next dense block. Clearly as the number of layers grow, it can become prohibitively expensive for information to propagate through the network. DenseNet can potentially overcome the vanishing gradient issue, but it is only an adhoc method \cite{huang2017denseNet, Zhang_2018}. Some other networks that have attempted to induce multilayer connections are Highway Net \cite{Srivastava2015HighwayNet}, AdaNet \cite{cortes2016adanet}, 
ResNetPlus \cite{ResNetPlus2018}, etc. All these models, however, largely lack rigorous mathematical frameworks. Furthermore, rigorous approaches to learn nonsmooth functions such as the absolute value function $|x|$ are scarce \cite{Imaizumi2018_nonsmooth}.

There has been a recent push in the scientific community to develop rigorous mathematical models and understanding of the DNNs \cite{Weinan_2019}. One way of doing so is to look at their architecture as dynamical systems. The articles \cite{haber2017stable, lu2017finite, ruthotto2018deep, Carola_resnet_2019, benning2019deep} have established that a DNN can be regarded as an optimization problem subject to a discrete ordinary differential equation (ODE) as constraints. The limiting problem in the continuous setting is an ODE constrained optimization problem \cite{ruthotto2018deep, Carola_resnet_2019}. Notice that  designing the solution algorithms at the continuous level can lead to architecture independence, i.e., the number of iterations remains the same even if the number of layers is increased.

The purpose of this paper is to present a novel fractional deep neural network 
which allows the network to access historic information of input and gradients across all subsequent layers. This is facilitated via our proposed use of fractional derivative based ODE as constraints.  We derive the optimality conditions for this network using the Lagrangian approach. Next, we consider a discretization for this fractional ODE and the resulting DNN is called \emph{Fractional-DNN}. We provide the algorithm and show numerical examples on some standard datasets.

Owing to the fact that fractional time derivatives allow memory effects, in the Fractional-DNN all the layers are connected to one another, with an appropriate scaling. In addition, fractional time  derivatives can be applied to nonsmooth functions \cite{antil_frac_time}. Thus, we aim to keep the benefits of standard DNN and the ideology of DenseNet, but remove the bottlenecks. 

The learning rate in a neural network is an important hyper-parameter which influences training \cite{Bengio2012LearningRate}. In our numerical experiments, we have observed an improvement in the learning rate via Fractional-DNN, which enhances the training capability of the network. Our numerical examples illustrate that, Fractional-DNN can potentially solve the vanishing gradient issue (due to memory), and handle nonsmooth data. 

The paper is organized as follows. In \cref{Prelim} we introduce notations and definitions. We introduce our proposed Fractional-DNN in \cref{fracDNN}. This is followed by \cref{NumericalApprox} where we discuss its numerical approximation. In \cref{Frac_DNN_alg}, we state our algorithm. The numerical examples given in \cref{NumExp} show the working and improvements due to the proposed ideas on three different datasets. 

\section{Preliminaries\label{Prelim}}

The purpose of this section is to introduce some notations and definitions that we will use throughout the paper. We begin with \cref{notations} where we state the standard notations. In \cref{s:ce} we describe
the well-known softmax loss function. \Cref{RL_deriv} is dedicated to the Caputo fractional time 
derivative. 
\begin{table}[!htb]
\centering
\caption{Table of Notations.\label{notations}}
\begin{tabular}{|c|l|}
\hline
\textbf{Symbol}			& \textbf{Description}\\ \hline
$n \in \NN$		& Number of distinct samples\\
$n_f \in \NN$  	& Number of sample features\\
$n_c \in \NN$  	& Number of classes\\
$N \in \NN$		& Number of network layers (i.e. network depth) \\
$Y \in \RR^{n_f \times n} $ & $Y = \{y^{(i)} \}_{i=1}^{n}$ is the collective feature set of $n$ samples. \\
$C_{obs} \in \RR^{n_c \times n}$ & $C_{obs} = \{c^{(i)}\}_{i=1}^{n}$ are the true class labels of the input data \\
$W \in \RR^{n_c \times n_f}$ & Weights \\
$K \in \RR^{n_f\; \times\; n_f}$ & Linear operator (distinct for each layer) \\
$b \in \RR$					 & Bias (distinct for each layer) \\
$P \in \RR^{n_f \times n} $ & Lagrange multiplier \\
$e_{n_c} \in \RR^{n_c}$		 & A vector of ones \\
$\tau \in \RR$								 & Time step-length\\
$\sigma(\cdot)$							 & Activation function, acting pointwise\\
$\gamma $					 & Order of fractional time derivative \\
$(\cdot)^{\prime}$					 & Derivative w.r.t. the argument \\
$tr(\cdot)$						 & Trace operator \\
$(\cdot)^{\intercal}$							 & Matrix transpose\\
$\odot$							 & Point-wise multiplication \\
$m_{1}$				 & Max count for randomly selecting a mini-batch in training 
\\
$m_{2}$				 & Max iteration count for gradient-based optimization solver 
\\
$\alpha_{train}, \; \alpha_{test}$ & Percentage of training and testing data correctly identified \\
\hline
\end{tabular}
\end{table}

\subsection{Cross Entropy with Softmax Function}
\label{s:ce}

Given collective feature matrix $Y$ 
with true labels  $C_{obs}$ and the unknown weights $W$, the cross entropy loss function given by
\begin{equation} \label{CE}
E(W,Y,C_{obs}) = -\frac{1}{n}\: \text{tr} (C^{\intercal}_{obs} \:\log(S(W,Y))) 
\end{equation}
measures the discrepancy between the true labels $C_{obs}$ and the predicted labels $\log(S(W,Y))$. Here, 
\begin{equation}\label{softmax}
S(W,Y) := \exp(WY)\:\text{diag}\left(\frac{1}{e^{\intercal}_{n_c} \exp(WY)} \right) 
\end{equation}
is the softmax classifier function, which gives normalized probabilities of samples belonging to the classes.

\subsection{Caputo Fractional Derivative} \label{RL_deriv}

In this section, we define the notion of Caputo fractional derivative and refer \cite{antil_frac_time} and references therein for the following definitions. 
	\begin{definition}[Left Caputo Fractional Derivative] For a fixed real number $0< \gamma < 1$, and an absolutely continuous function $u\colon [0,T]\rightarrow \mathbb{R}$, the left Caputo fractional derivative is defined by:
	\begin{equation} \label{Lcaputo}
	d^{\gamma}_{t} u(t) 
	=\frac{1}{\Gamma(1-\gamma)} \frac{d}{dt}\int_{0}^{t}\frac{u(r)-u(0)}{(t-r)^{\gamma}}\;dr,
	\end{equation}
	where $\Gamma (\cdot)$ is the Euler-Gamma function.
	\end{definition}
	\begin{definition}[Right Caputo Fractional Derivative] For a fixed real number $0< \gamma < 1$, and an absolutely continuous function $u\colon [0,T]\rightarrow \mathbb{R}$, the right Caputo fractional derivative is defined by:
		\begin{equation} \label{Rcaputo}
		d^{\gamma}_{T-t} u(t) =
		\frac{-1}{\Gamma(1-\gamma)} \frac{d}{dt}\int_{t}^{T}\frac{u(r)-u(T)}{(r-t)^{\gamma}}\;dr.
		\end{equation}
	\end{definition}
Notice that, $d^{\gamma}_{t} u(t)$ and $d^{\gamma}_{T-t} u(t)$ in definitions \cref{Lcaputo} and \cref{Rcaputo} exist almost everywhere on $[0,T]$, \cite[Theorem 2.1]{kilbas_srivastava_fde}, and are represented, respectively, by
\[ d^{\gamma}_{t} u(t) 
=\frac{1}{\Gamma(1-\gamma)} \int_{0}^{t}\frac{u'(r)}{(t-r)^{\gamma}}\;dr,\qquad \mbox{and} \qquad  d^{\gamma}_{T-t} u(t) 
=\frac{-1}{\Gamma(1-\gamma)} \int_{t}^{T}\frac{u'(r)}{(r-t)^{\gamma}}\;dr. \]
Moreover, if $\gamma=1$ and $u \in C^1([0,T])$, then one can show that $d^{\gamma}_{t} u(t) = u'(t) = d^{\gamma}_{T-t} u(t)$. 
We note that the fractional derivatives in \cref{Lcaputo} and \cref{Rcaputo} are nonlocal operators. Indeed, the derivative of $u$ at a point $t$ depends on all the past and future events, respectively. This behavior is different than the classical case of $\gamma=1$. 

The left and right Caputo fractional derivatives are linked by the fractional integration by parts formula, ~\cite[Lemma 3]{antil_otarola_salgado}, which will be stated next.
For $\gamma \in (0,1)$, let
\[ 
	\mathbb{L}_{\gamma}\coloneqq \left\{ f\in C([0,T]): d_t^{\gamma}f\in L^2(0,T) \right\}, \qquad \mathbb{R}_{\gamma}\coloneqq \left\{ f\in C([0,T]): d_{T-t}^{\gamma}f\in L^2(0,T) \right\} . 
\]
\begin{lemma}[Fractional Integration-by-Parts] \label{int_by_parts}
	For $f\in \mathbb{L}_{\gamma}$ and $g\in \mathbb{R}_{\gamma}$, the following integration-by-parts formula holds:	
	\begin{equation}\label{frac_IP}
		\int_{0}^{T} d^{\gamma}_{t} f(t) g(t) \; dt= \int_{0}^{T}  f(t) d^{\gamma}_{T-t} g(t) \; dt + g(T)(I_t^{1-\gamma}f)(T) - f(0) (I_{T-t}^{1-\gamma}g)(0) , 
	\end{equation}
	where $I_t^{1-\gamma}w(t)$ and $I_{T-t}^{1-\gamma}w(t)$ are the left and right Riemann-Liouville fractional integrals of order $\gamma$ and are given by 
	\[  
		I_t^{1-\gamma}w(t):=\frac{1}{\Gamma(1-\gamma)} \int_{0}^{t}\frac{w(r)}{(t-r)^{\gamma}}\;dr \; \mbox{ and } \; I_{T-t}^{1-\gamma}
		w(t):=\frac{1}{\Gamma(1-\gamma)} \int_{t}^{T}\frac{w(r)}{(r-t)^{\gamma}}\;dr. 
	\]
\end{lemma}
\section{Continuous Fractional Deep Neural Network\label{fracDNN}}

After the above preparations, in this section, we shall introduce the Fractional-DNN. First we briefly describe the classical
RNN, and then extend it to develop the Fractional-DNN. We formulate  our problem as a constrained optimization problem. Subsequently, 
we shall use the Lagrangian approach to derive the optimality conditions.
\subsection{Classical RNN}

Our goal is to approximate a map $\mathcal{F}$.
A classical RNN helps approximate $\mathcal{F}$, for a known set of inputs and outputs. To construct an RNN, 
for each layer $j$, we first consider a linear-transformation of $Y_{j-1}$ as, 
\[
\mathcal{G}_{j-1}(Y_{j-1}) = K_{j-1}Y_{j-1} + b_{j-1},
\]
where the pair $(K_{j},b_{j})$ denotes an unknown linear operator and bias at the $j^{th}$ layer. When $N>1$ then the network is considered ``deep". Next
we introduce non linearity using a nonlinear activation function $\sigma$ (e.g. ReLU or $\tanh$). 
The resulting RNN is,
\begin{equation}\label{eq:genitdisc}
Y_j = Y_{j-1} + \tau (\sigma \circ \mathcal{G}_{j-1})(Y_{j-1}), \quad j = 1,\cdots,N; \; \;\;  N>1,
\end{equation}
where $\tau>0$ is the time-step. 
Finally, the RNN approximation of $\mathcal{F}$ is given by,
\[
\mathcal{F}_\theta(\cdot) = \Big( 
\big(I +\tau (\sigma \circ \mathcal{G}_{{N-1}})\big) \circ \big(I +\tau (\sigma \circ \mathcal{G}_{{N-2}})\big) \circ \dots \circ \big(I +\tau (\sigma \circ \mathcal{G}_{0})\big)
\Big)(\cdot),
\]
with $\theta = (K_j,b_j)$ as the unknown parameters. In other words, the problem of approximating $\mathcal{F}$ using classical RNN, intrinsically, is a problem of learning $(K_j,b_j)$. 

Hence, for given datum $(Y_0,C)$, the learning problem then reduces to minimizing a loss function $J(\theta,(Y_N,C))$, 
subject to constraint \cref{eq:genitdisc}, i.e., 
\begin{equation}\label{eq:discNN}
\begin{aligned}
&\min_{\theta} \;\; \mathcal{J}(\theta,(Y_N,C))\\
\text{s.t.}\quad Y_j &= Y_{j-1} + \tau (\sigma \circ \mathcal{G}_{j-1})(Y_{j-1}), \quad j = 1,\dots,N.
\end{aligned}
\end{equation}
Notice that the system \cref{eq:genitdisc} is the forward-Euler discretization of the following continuous in time ODE, 
see \cite{he2016deep, gunther2018layer, ruthotto2018deep}, 
\begin{equation}\label{eq:genitcont}
\begin{aligned}
 d_t Y(t)    &= \sigma (K(t) Y(t) + b(t) )  , \quad t \in (0,T),  \\ 
	Y(0) &= Y_0 .
\end{aligned}	
\end{equation}	
The  continuous learning problem then requires minimizing the loss function $\mathcal{J}$ at the final time $T$
subject to the ODE constraints \cref{eq:genitcont}:
\begin{equation} \label{cont_gen_prob}
\begin{aligned}
\min_{\theta = (K,b)} \;\;\; &\mathcal{J}(\theta,(Y(T),C))\\
\text{s.t.}\quad 	& \cref{eq:genitcont}  
\end{aligned}
\end{equation}

Notice that designing algorithms for the continuous in time problem \cref{cont_gen_prob} instead of the discrete in time problem \cref{eq:discNN} has several key advantages. In particular, it will  
lead to algorithms which are independent of the neural network architecture, i.e., independent of the number of 
layers. In addition, the approach of \cref{cont_gen_prob} can help us determine the stability of the neural network \cref{eq:discNN}, see \cite{benning2019deep, haber2017stable}. 
Moreover, for the neural network \cref{eq:discNN}, it has been noted that as the information about the input or gradient passes through many layers, it can vanish and ``wash out'', or grow and ``explode'' exponentially \cite{Bengio1994}. 
There have been adhoc attempts to address these concerns, see for instance \cite{Srivastava2015HighwayNet, cortes2016adanet, huang2017denseNet}, but a satisfactory mathematical explanation and model does not currently exist. One of the main goals of this paper is to introduce such a model.

Notice that \cref{eq:genitcont}, and its discrete version \cref{eq:genitdisc}, incorporates many algorithmic processes such as linear solvers, preconditioners, nonlinear solvers, optimization solvers, etc. Furthermore, there are well-established numerical algorithms that re-use information from previous iterations to accelerate convergence, e.g. the BFGS method \cite{NocedalWright2006}, Anderson acceleration \cite{anderson1965iterative}, and variance reduction methods \cite{roux2012stochastic}. These methods account for the history $Y_j, Y_{j-1}, Y_{j-2}, \dots, Y_0$, while choosing $Y_{j+1}$. Motivated by these observations we  introduce versions of \cref{eq:genitdisc} and \cref{eq:genitcont} that can account for history (or memory) effects in a rigorous mathematical fashion. 	

\subsection{Continuous Fractional-DNN}

The fractional time derivative in \cref{Lcaputo} has a distinct ability to allow a memory effect, for instance in materials with hereditary properties \cite{brown2018analysis}. Fractional time derivative can be derived by using the anomalous random walks where the walker experiences delays between jumps \cite{metzler2000random}. In contrast, 
the standard time derivative naturally arises in the case of classical random walks. We use the idea of fractional time derivative  
to enrich the constraint optimization problem  \cref{cont_gen_prob}, and subsequently \cref{eq:discNN}, by replacing the standard time derivative $d_t$ by the fractional time derivative $d_t^\gamma$ of order $\gamma \in (0,1)$. Recall that for $\gamma = 1$, we obtain the classical derivative $d_t$. 
Our new continuous in time model, the Fractional-DNN, is then given by (cf. \cref{eq:genitcont}),
\begin{equation}\label{eq:b}
\begin{aligned}
d_t^\gamma Y(t) &= \mathcal{F}_{\theta}(Y(t),t,\theta(t)), \quad t \in (0,T),  \\
Y(0) &= Y_0 
\end{aligned}
\end{equation}
where $d_t^\gamma$ is the Caputo fractional derivative as defined in \cref{Lcaputo}. 
The discrete formulation 
of Fractional-DNN will be discussed in the subsequent section. 

The main reason for using the Caputo fractional time derivative over its other counterparts such as the Riemann Liouville fractional derivative is the fact that the Caputo derivative of a constant function is zero and one can impose the initial conditions $Y(0) = Y_0$ in a classical manner \cite{Samko1993}. Note 
that $d_t^\gamma$ is a nonlocal operator in a sense that in order to evaluate the fractional derivative of $Y$ at a point $t$, we need the cumulative
information of $Y$ over the entire sub-interval $[0,t)$. This is how the Fractional-DNN enables connectivity across all antecedent layers (hence the memory effect). As we shall illustrate with the help of a numerical example in \cref{NumExp}, this feature can help overcome the vanishing gradient issue, as the cumulative effect of the gradient of the precedent layers is less likely to be zero. 

\begin{remark}[Caputo Derivative of Nonsmooth Functions]
\label{rem:nonsmooth_deriv}
{\rm 
The Caputo fractional derivative \cref{Lcaputo} can be applied to non-smooth functions. 
Consider, e.g.
$Y(t) \coloneqq |t|.$
Notice that $Y(t)$ is not differentiable at $t = 0$. However, \cref{Lcaputo} 
yields, 
$d^{\gamma}_{t} Y(t) = \frac{1}{\Gamma(2-\gamma)} t^{1-\gamma}. $
Since $\gamma \in (0,1)$, therefore $d^{\gamma}_{t} Y(t)$ at $t = 0$ is zero.
}
\qed 
\end{remark}

Owing to \cref{rem:nonsmooth_deriv} we can better account for features, $Y$, 
which are non-smooth, as a result of which the smoothness requirement on the unknown parameters $\theta$ can be weakened. 
This, in essence, can help with the exploding gradient issue in DNNs. 
  
The generic learning problem with Fractional-DNN as constraints can be expressed as,
\begin{equation} \label{cont_gen_prob_frac}
\begin{aligned}
\min_{\theta = (K,b)}  \;\;\; &\mathcal{J}(\theta,(Y(T),C))\\
\text{s.t.}\quad 	& \cref{eq:b} 
\end{aligned}
\end{equation}
Note that the choice of $\mathcal{J}$ depends on the type of learning problem. We will next consider a specific structure 
of $\mathcal{J}$ given by the cross entropy loss functional, defined in \cref{CE}.

\subsection{Continuous Fractional-DNN and Cross Entropy Loss Functional}
Supervised learning problems are a broad class of machine learning problems which use labeled data. These problems are further divided into two types, namely regression problems and classification problems. The specific type of the problem dictates the choice of $\mathcal{J}$ in \cref{cont_gen_prob_frac}. Regression problems often occur in physics informed models, e.g. sample reconstruction inverse problems \cite{ADK_19_ml_tomo, Pock2017_MRI}. On the other hand, classification problems occur, for instance, in computer vision \cite{scherer2020, CIRESAN2012333}. In both the cases, a neural network is used to learn the unknown parameters.  In the discussion below we shall focus on classification problems, however, the entire discussion directly applies to 
regression type problems.

Recall that the cross entropy loss functional $E$, defined in \cref{CE}, measures the discrepancy between the actual and the predicated classes. Replacing, $\mathcal{J}$ in \cref{cont_gen_prob_frac} by $E$ together with a regularization term 
$\mathcal{R}(W,K(t),b(t))$, we arrive at
\begin{equation} \label{RNN_frac}
\begin{aligned}
\min_{W,K,b} \hspace{0.5cm} &E(W,Y(T),C_{obs}) + \R(W,K(t),b(t)) \\
\text{s.t.} \hspace{0.5cm}& \left\{
\begin{aligned}
d^{\gamma}_t Y(t) &= \: \sigma(K(t)Y(t)+b(t)), \hspace{1cm} t\in (0,T), \\
Y(0) &= \: Y_0 \, .
\end{aligned}
\right.
\end{aligned}
\end{equation}
Note that, in this case, the unknown parameter $\theta \coloneqq(W,K,b)$, where $K$ and $b$ are, respectively, the linear operator and bias for each layer, and the weights $W$ are a feature-to-class map. Furthermore, $\sigma$ is a nonlinear activation function and $(Y_0,C_{obs})$ is the given data, with $C_{obs}$ as the true labels of $Y_0$. 

To solve \cref{RNN_frac}, we rewrite this problem as an unconstrained optimization problem via the Lagrangian functional and derive the optimality conditions. Let $P$ denote the Lagrange multiplier, then the Lagrangian functional is given by, 
\begin{equation*}
\mathcal{L}(Y,W,K,b;P) := \;E(W,Y(T),C_{obs}) + \R(W,K(t),b(t)) + \langle d^{\gamma}_t Y(t) - \sigma(K(t)Y(t)+b(t)),P(t) \rangle  ,
\end{equation*}
where, $\langle \cdot, \cdot \rangle := \int_0^T \langle \cdot, \cdot \rangle_F \;dt$ is the $L^2$-inner product, and $\langle \cdot, \cdot \rangle_F$ is the Frobenius 
inner product. Using the fractional integration-by-parts from \cref{frac_IP}, we obtain
\begin{equation} \label{Lag}
\begin{aligned} 
\mathcal{L}(Y,W,K,b;P) 
=\: &E(W,Y(T),C_{obs}) + \R(W,K(t),b(t)) -\langle \sigma(K(t)Y(t)+b(t)),P(t)\rangle \\
& + \langle  Y(t), d^{\gamma}_{T-t} P(t)\rangle + \langle (I_t^{1-\gamma}Y)(T),P(T) \rangle_F - \langle Y_0, (I_{T-t}^{1-\gamma}P)(0)\rangle_F.
\end{aligned}
\end{equation}  

Let $(\overline{Y},\overline{W},\overline{K},\overline{b};\overline{P})$ denote a stationary point, then the first order necessary optimality conditions are given by the following set of 
state, adjoint and design equations:
\begin{enumerate}[label=(\Alph*)]
\item \textbf{State Equation.} The gradient of $\mathcal{L}$ with respect to $P$ at $(\overline{Y},\overline{W},\overline{K},\overline{b};\overline{P})$ yields the state equation $\nabla_{P} \mathcal{L}(\overline{Y},\overline{W},\overline{K},\overline{b};\overline{P})= 0$, equivalently, 
\begin{equation}\label{eq:state_cont}
\begin{aligned}
d^{\gamma}_t \overline{Y}(t) &= \: \sigma(\overline{K}(t)\overline{Y}(t)+\overline{b}(t)), \hspace{1cm} t\in (0,T),  \\
\overline{Y}(0) &= \: Y_0 
\end{aligned}
\end{equation}
where $d^{\gamma}_t$ denotes the left Caputo fractional derivative \cref{Lcaputo}. 
In \cref{eq:state_cont}, for the state variable $\overline{Y}$, we solve forward in time, therefore we call 
\cref{eq:state_cont}  as the \textit{forward propagation}.

\item \textbf{Adjoint Equation.} Next, the gradient of $\mathcal{L}$ with respect to $Y$ at $(\overline{Y},\overline{W},\overline{K},\overline{b};\overline{P})$ yields the adjoint equation $\nabla_{Y} \mathcal{L}(\overline{Y},\overline{W},\overline{K},\overline{b};\overline{P})= 0$, equivalently, 
\begin{equation}\label{eq:adjnt_cont}
\begin{aligned}
d^{\gamma}_{T-t}{\overline{P}}(t) &= (\sigma^{\prime} (\overline{K}(t)\overline{Y}(t)+\overline{b}(t))\; \overline{K}(t))^{\intercal} \; \overline{P}(t)  \\
&=\overline{K}(t)^{\intercal}\left(\overline{P}(t) \odot \sigma^{\prime}\left(\overline{K}(t) \overline{Y}(t) + \overline{b}(t)\right)\right),\hspace{1cm} t\in (0,T),\\
\overline{P}(T) & 
= -\frac{1}{n}\overline{W}^{\intercal}(-C_{obs}+ S(\overline{W},\overline{Y}(T)))
\end{aligned}
\end{equation}
where $d^{\gamma}_{T-t}$ denotes the right Caputo fractional derivative \cref{Rcaputo} 
and $S$ is the softmax function defined in \cref{softmax}. Notice that the adjoint variable $\overline{P}$ in  \cref{eq:adjnt_cont}, with its terminal condition, is obtained by marching backward in time. As a result, the equation \cref{eq:adjnt_cont} is called \textit{backward propagation}. 

\item \textbf{Design Equations.} Finally, equating $\nabla_{W}  \mathcal{L}(\overline{Y},\overline{W},\overline{K},\overline{b};\overline{P})$, $\nabla_{K}  \mathcal{L}(\overline{Y},\overline{W},\overline{K},\overline{b};\overline{P})$, and $\nabla_{b}  \mathcal{L}(\overline{Y},\overline{W},\overline{K},\overline{b};\overline{P})$ to zero, respectively, yields the design equations (with 

($\overline{W}, \overline{K}, \overline{b}$) as the design variables),

\begin{alignat}{4} \label{eq:design_cont}	
	\nabla_{W}  \mathcal{L}(\overline{Y},\overline{W},\overline{K},\overline{b};\overline{P}) 
	= &\frac{1}{n}\big(-C_{obs}+ S(\overline{W},\overline{Y}(T))\big)\;\big(\overline{Y}(T)\big)^{\intercal}  \nonumber \\ 
	&+ \nabla_W \R(\overline{W},\overline{K}(T),\overline{b}(T)) =  0, \nonumber\\	
	\nabla_{K}  \mathcal{L}(\overline{Y},\overline{W},\overline{K},\overline{b};\overline{P})= &- \overline{Y}(t)\:\left(\overline{P}(t) \odot \sigma^{\prime}(\overline{K}(t) \overline{Y}(t) + \overline{b}(t))\right)^{\intercal}  \\ 
	& + \nabla_K \R(\overline{W},\overline{K}(t),\overline{b}(t))  = 0, \nonumber\\
	\nabla_{b}  \mathcal{L}(\overline{Y},\overline{W},\overline{K},\overline{b};\overline{P})= &- \langle \sigma'(\overline{K}(t)\overline{Y}(t)+\overline{b}(t)), \overline{P}(t)\rangle_F  \nonumber \\ 
	&+  \nabla_{b} \R(\overline{W},\overline{K}(t),\overline{b}(t))  = 0, \nonumber
\end{alignat}
for almost every $t \in (0,T)$. 
\end{enumerate}

In view of (A)-(C), we can use a gradient based solver to find a stationary point to \cref{RNN_frac}.
\begin{remark}{\rm(Parametric Kernel $K(\psi(t))$).
Throughout our discussion, we have assumed $K(t)$ to be some unknown linear operator. We remark that a structure could also be prescribed to $K(t)$, parameterized by a stencil $\psi$. Then, the kernel is $K(\psi(t))$, and the design variables now are $\theta = (W,\psi,b)$. Consequently, $K(\psi(t))$ can be thought of as a differential operator on the feature space, e.g. discrete Laplacian with a five point stencil. It then remains to compute the sensitivity of the Lagrangian functional w.r.t. $\psi$ to get the design equation.  Note that this approach can further reduce the number of unknowns. \qed 
}
\end{remark}
Notice that so far the entire discussion has been at the continuous level and it has been independent of the number of 
network layers. Thus, it is expected that if we discretize (in time) the above optimality system, then the resulting gradient 
based solver is independent of the number of layers. We shall discretize the above optimality system in the next section. 


\section{Discrete Fractional Deep Neural Network\label{NumericalApprox}}
 \hspace{-2pt}We shall adopt the \emph{optimize-then-discretize} approach. Recall that the first order stationarity conditions for the 
continuous problem \cref{RNN_frac} are given in \cref{eq:state_cont}, \cref{eq:adjnt_cont}, and \cref{eq:design_cont}.
In order to discretize this system of equations, we shall first discuss the approximation of 
Caputo fractional derivative. 

\subsection{Approximation of Caputo Derivative} \label{RLGL}
There exist various approaches to discretize the fractional Caputo derivative. We will use the $L^1$-scheme 
\cite{antil_otarola_salgado, Mstynes} to discretize the left and right Caputo fractional derivative 
$d_t^{\gamma}u(t)$ and $d_{T-t}^{\gamma}u(t)$ given in \cref{Lcaputo} and \cref{Rcaputo}, respectively.

Consider the following fractional differential equation involving the \textbf{left} Caputo fractional derivative, for $0< \gamma <1$,
\begin{equation} \label{eq:eq1}
d_t^{\gamma}u(t)=f(u(t)), \quad u(0)=u_0.
\end{equation}
We begin by discretizing the time interval $[0,T]$ uniformly with step size $\tau$,
\[ 
0=t_0<t_1<t_2<\cdots <t_{j+1}<\cdots< t_N = T, \;\;\mbox{where }\; t_j=j \tau .
\]	
Then using the $L^1$-scheme, the discretization of \cref{eq:eq1} is given by
\begin{alignat}{3} \label{eq:disc_cCaputo}
u(t_{j+1})= u(t_{j}) - \sum_{k=0}^{j-1}   a_{j-k}\: \left(u(t_{k+1}) - u(t_k) \right)  + \tau^{\gamma}\Gamma(2-\gamma)f(u(t_{j})).
\quad  j = 0,...,N-1 \;,
\end{alignat}
where coefficients $a_k$ are given by, 
\begin{equation}\label{a_k}
a_{j-k}=(j+1-k)^{1-\gamma}-(j-k)^{1-\gamma}.
\end{equation}

Next, let us consider the discretization of the fractional differential equation involving the \textbf{right} Caputo fractional operator, for $0< \gamma <1$,
\begin{equation} \label{eq:eq2}
d_{T-t}^{\gamma}u(t)=f(u(t)), \quad u(T)=u_T. 
\end{equation}
Again using $L^1$-scheme we get the following discretization of \cref{eq:eq2}:
\begin{alignat}{3} \label{eq:disc_cRCaputo}
u(t_{j-1})= u(t_{j}) + \sum_{k=j}^{N-1}   a_{k-j}\: \left(u(t_{k+1}) - u(t_k) \right)  - \tau^{\gamma}\Gamma(2-\gamma)f(u(t_j)).
\quad  j = N,...,1.
\end{alignat}

The example below illustrates a numerical implementation of the $L^1$-scheme \cref{eq:disc_cCaputo}.
\begin{example}
Consider the linear differential equation
	\begin{alignat}{2} \label{ex:mittag}
		d_t^{\,0.5} u(t)=-4 u(t),\;\; u(0)=0.5.
	\end{alignat} 
	Then, the solution to \cref{ex:mittag} is given by, see \cite[Section 42]{Samko1993}, also \cite[Section 1.2]{Podlubny1999}
	\begin{alignat}{3}
			u(t)=0.5\: E_{0.5}(-4t^{0.5}),
	\end{alignat}
	where $E_{\alpha}$, with $\alpha >0$, is the Mittag Leffler function defined by 
	\[ E_{\alpha }(z)=E_{\alpha , 1}(z)= \sum_0^{\infty} \frac{z^k}{\Gamma(\alpha k +1)} .\]
	\Cref{fig:GL_approx} depicts the true solution and the numerical solutions using discretization \cref{eq:disc_cCaputo} for the above example with uniform step size $\tau=0.005$ and final time, $T=1$. 
	\begin{figure}[htb] 
		\begin{center}
			\includegraphics[width=0.5\textwidth]{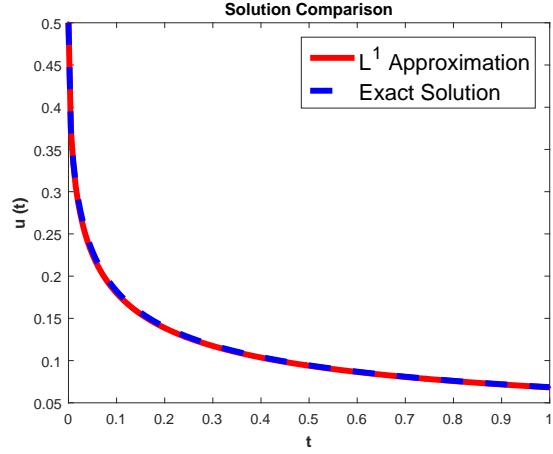}
			\caption{Comparison of the exact solution of \cref{ex:mittag} (\textit{blue}) with an $L^1$ scheme approximation (\textit{red}).}
			\label{fig:GL_approx} 
		\end{center}
	\end{figure}
\qed
\end{example}

\subsection{Discrete Optimality Conditions}

Next, we shall discretize the optimality conditions given in \cref{eq:state_cont} -- \cref{eq:design_cont}.
Notice that, each time-step corresponds to one layer of the neural network. It is necessary to do one forward
propagation (state solve) and one backward propagation (adjoint solve) to derive an expression of the 
gradient with respect to the design variables.

\begin{enumerate} [label=(\Alph*)]
\item \textbf{Discrete State Equation.} We use the $L^1$ scheme discussed in \cref{eq:disc_cCaputo} to 
discretize the state equation \cref{eq:state_cont} and arrive at 
\begin{equation}\label{fwd_prop_L1}
	\begin{aligned}
	\overline{Y}(t_{j}) \;=\; &Y(t_{j-1}) - \sum_{k=1}^{j-1}   a_{j-k}\: \left(Y(t_{k}) - Y(t_{k-1}) \right)  \\
	&\hspace{1.5cm} +  \tau^{\gamma}\Gamma(2-\gamma)\sigma(\overline{K}(t_{j-1}) \overline{Y}(t_{j-1}) + \overline{b}(t_{j-1})), \qquad j = 1,...,N\\
	\overline{Y}(t_0) \;=\; &Y_0
	\end{aligned}
\end{equation}

\item \textbf{Discrete Adjoint Equation.} We use the $L^1$ scheme discussed in \cref{eq:disc_cRCaputo} to 
discretize the adjoint equation \cref{eq:adjnt_cont} and arrive at
\begin{equation}\label{bckwd_prop_L1}
\begin{aligned}
\overline{P}(t_{j}) &= P(t_{j+1}) + \sum_{k=j+1}^{N-1} a_{k-j-1}\: \left(P(t_{k+1}) - P(t_k) \right)- \hspace{2cm} j = N-1,...,0\\
&\hspace{1.5cm} \tau^{\gamma}\Gamma(2-\gamma)\left[-\overline{K}(t_{j})^{\intercal}\left(\overline{P}(t_{j+1}) \odot \sigma^{\prime}\left(\overline{K}(t_{j}) \overline{Y}(t_{j+1}) + \overline{b}(t_{j})\right)\right)\right], \\
\overline{P}(t_N) 
&= -\frac{1}{n}\overline{W}^{\intercal}(-C_{obs}+ S(\overline{W},\overline{Y}(t_N))) 
\end{aligned}
\end{equation}

\item \textbf{Discrete Gradient w.r.t. Design Variables.} For $j=0,\dots,N-1$, the approximation of the gradient 
\cref{eq:design_cont} with respect to the design variables is given by,
\begin{equation}\label{eq:design_disc}
\begin{aligned}
	\nabla_W \mathcal{L}(\overline{Y},\overline{W},\overline{K},\overline{b};\overline{P}) = &\frac{1}{n}\big(-C_{obs}+ S(\overline{W},\overline{Y}(t_N))\big)\;\big(\overline{Y}(t_N)\big)^{\intercal}  \\
 &+\nabla_W \R(\overline{W},\overline{K}(t_N),\overline{b}(t_N))\\
\nabla_{K} \mathcal{L}(\overline{Y},\overline{W},\overline{K},\overline{b};\overline{P}) = &-\overline{Y}(t_j)\:\left(\overline{P}(t_{j+1}) \odot \sigma^{\prime}(\overline{K}(t_j) \overline{Y}(t_j) + \overline{b}(t_j))\right)^{\intercal}  \\
 &+\nabla_K \R(\overline{W},\overline{K}(t_j),\overline{b}(t_j))\\
\nabla_{b} \mathcal{L}(\overline{Y},\overline{W},\overline{K},\overline{b};\overline{P}) = &-\; \langle \sigma^{\prime}(\overline{K}(t_j) \overline{Y}(t_j) + \overline{b}(t_j)),\overline{P}(t_{j+1})\rangle_F  \\
 &+ \nabla_b \R(\overline{W},\overline{K}(t_j),\overline{b}(t_j)) \; .
\end{aligned}
\end{equation}
\end{enumerate}

Whence, we shall create a gradient based method to solve the optimality condition \cref{fwd_prop_L1}-\cref{eq:design_disc}. We reiterate that each computation of the gradient in \cref{eq:design_disc}, requires one state and one adjoint solve.

\section{Fractional-DNN Algorithm\label{Frac_DNN_alg}}

Fractional-DNN is a supervised learning architecture, i.e. it comprises of a training phase and a testing phase. During the training phase, 
labeled data is passed into the network and the unknown parameters are learnt. Those parameters then define the trained Fractional-DNN model for that type of data. Next, a testing dataset, which comprises of data previously unseen by the network, is passed to the trained net, and a prediction of classification is obtained. This stage is known as the testing phase. Here the true classification is not shown to the network when a prediction is being made, but can later be used to compare the network efficiency, as we have done in our numerics. The three important components of the algorithmic structure are forward propagation, backward propagation, and gradient update. The forward and backward propagation structures are given in \cref{alg:fDNN_fwd_prop_L1,alg:fDNN_bckwd_prop_L1}. The gradient update is accomplished in the training phase, discussed in \cref{sec:train_phase}. Lastly, the testing phase of the algorithm is discussed in \cref{sec:test_phase}. 
%
\begin{algorithm}
\caption{Forward Propagation in Factional-DNN ($L^1$-scheme)}
\label{alg:fDNN_fwd_prop_L1}
\begin{algorithmic}[1]
\REQUIRE $\left(Y_0,C_{obs}\right), W,\: \{K_j, b_j\}_{j=0}^{N-1}, \:N, \:\tau, \:\gamma$ 
\ENSURE $\{Y_j\}_{j=1}^{N},\; P_N$, 
\STATE Let  $z_0 = 0$. 
\FOR {$j = 1,\cdots,N$}
\FOR {$k = 1,\cdots,j-1$} 
\STATE Compute $a_{j-k}$: \quad 
\COMMENT {Use \cref{a_k}}
\STATE Update $z_{k}$: \hspace{1cm} 
$
z_{k} = z_{k-1} + a_{j-k}\;(Y_{k}-Y_{k-1})
$
\ENDFOR 
\STATE Update $Y_{j}$: \hspace{1cm}
$
Y_{j} \;=\; Y_{j-1} - z_{j-1} + (\tau)^{\gamma}\:\Gamma(2-\gamma)\: \sigma(K_{j-1} Y_{j-1} + b_{j-1})
$
\ENDFOR
\STATE Compute $P_N$: \hspace{1cm}
$
P_N = -(n)^{-1} \; W^{\intercal}(-C_{obs}+ S(W,Y_N))
$
\end{algorithmic}
\end{algorithm}

\begin{algorithm}
\caption{Backward Propagation in Factional-DNN ($L^1$-scheme)}
\label{alg:fDNN_bckwd_prop_L1}
\begin{algorithmic}[1]
\REQUIRE $\{Y_j\}_{j=1}^{N},\; P_N,\: \{K_j, b_j\}_{j=0}^{N-1}, \:N, \:\tau, \:\gamma$ 
\ENSURE $\{P_j\}_{j=0}^{N-1}$ 
\STATE Let $x_0 = 0$.
\FOR {$j = N-1,\cdots,0$} 
\FOR {{$k = j+1,\cdots,N-1$}}
\STATE Compute $a_{k-j-1}$: \quad
\COMMENT {Use \cref{a_k}}
\STATE Compute $x_{k}$: \quad
$
x_k = x_{k-1} + a_{k-j-1}\;(P_{k+1} - P_{k})
$
\ENDFOR
\STATE Update $P_j$: \quad$
P_{j} = P_{j+1}+ x_{N-1} - (\tau)^{\gamma}\:\Gamma(2-\gamma)\:[-K_{j}^{\intercal}(P_{j+1} \odot \sigma^{\prime}(K_{j} Y_{j+1} + b_{j}))]$
\ENDFOR
\end{algorithmic}
\end{algorithm}
\subsection{Training Phase}\label{sec:train_phase}

The training phase of Fractional-DNN is shown in  \cref{alg:fDNN_train}.

\begin{algorithm}[!htb]
\caption{Training Phase of Factional-DNN}
\label{alg:fDNN_train}
\begin{algorithmic}[1]
\REQUIRE $\left(Y_0,C_{obs}\right), \;N, \;\tau, \;\gamma,\;m_1,m_2$ 
\ENSURE $W,\: \{K_j, b_j\}_{j=0}^{N-1}, \; C_{train}, \alpha_{train}$, 
\STATE Initialize $W,\{K_j,b_j\}_{j=0}^{N-1}$ 
\FOR {$i = 1,\cdots,m_1$}
\STATE Let $(\hat{Y}_0,\hat{C}_{obs}) \subset \left(Y_0,C_{obs}\right)$ 
\COMMENT{Randomly select a mini-batch and apply BN using \cref{eq:BN}}
\STATE \textbf{FORWARD PROPAGATION} 
\COMMENT {Use \cref{alg:fDNN_fwd_prop_L1} to get $\{\hat{Y}_j\}_{j=1}^{N},\; P_N$}.
\STATE \textbf{BACKWARD PROPAGATION} 
\COMMENT {Use \cref{alg:fDNN_bckwd_prop_L1} to get $\{P_j\}_{j=0}^{N-1}}$.
\STATE \textbf{GRADIENT COMPUTATION}
\STATE Compute $\nabla_W \mathcal{L},\;\{\nabla_{K}\mathcal{L}\} ,\;\{\nabla_{b}\mathcal{L}\}$ 
\[
\begin{aligned}
\nabla_W \mathcal{L} \;=\; &(n)^{-1}\left(-C_{obs}+ S(W,\hat{Y}_N)\right)\;(\hat{Y}_N)^{\intercal} +\nabla_W \R(W,K_j,b_j)\\
\nabla_K \mathcal{L} \;=\; & -\hat{Y}_j\:\left(P_{j+1} \odot \sigma^{\prime}(K_j \hat{Y}_j + b_j)\right)^{\intercal} + \nabla_K \R(W,K_j,b_j)\\
\nabla_b \mathcal{L} \;=\; &-tr\left(\sigma^{\prime}(K_j \hat{Y}_j + b_j)\: P_{j+1}\right) + \nabla_b \R(W,K_j,b_j)
\end{aligned}
\]
\STATE Pass $\nabla_W \mathcal{L},\:\nabla_{K}\mathcal{L},\;\nabla_{b}\mathcal{L}$ to gradient based solver with $m_2$ max iterations to update $W, \{K_j,b_j\}_{j=0}^{N-1}$. 
\STATE Compute $\hat{C}_{train} = S(W,\hat{Y}_N)$
\STATE Compare $\hat{C}_{train}$ to $\hat{C}_{obs}$ to compute $\alpha_{train}$ 
\ENDFOR
\end{algorithmic}
\end{algorithm}
\subsection{Testing Phase}\label{sec:test_phase}
The testing phase of Fractional-DNN is shown in \cref{alg:fDNN_test}.

\begin{algorithm} [!htb]
\caption{Testing Phase of Fractional-DNN}
\label{alg:fDNN_test}
\begin{algorithmic}[1]
\REQUIRE $\left(Y_0^{test},C_{obs,test}\right),\; W,\: \{K_j, b_j\}_{j=0}^{N-1},\; N,\; \tau, \;\gamma$ 
\ENSURE  $C_{test}, \alpha_{test}$
\STATE Let  $Y_0 = Y_0^{test}$ 
\COMMENT {Apply BN using \cref{eq:BN}}
\STATE \textbf{FORWARD PROPAGATION} 
\COMMENT {Use \cref{alg:fDNN_fwd_prop_L1} to get $\{Y_j\}_{j=1}^{N}$}.
\STATE Compute $C_{test} = S(W,Y_N)$
\STATE Compare $C_{test}$ to $C_{obs,test}$ to compute $\alpha_{test}$ 
\end{algorithmic}
\end{algorithm}

\section{Numerical Experiments\label{NumExp}}

In this section, we present several numerical experiments where we use our proposed Fractional-DNN algorithm from \cref{Frac_DNN_alg} to solve classification problems for two different datasets. We recall that the goal of classification problems, as the name suggests, is to classify objects into pre-defined class labels. First we prepare a training dataset and along-with its classification, pass it to the training phase of Fractional-DNN (\cref{alg:fDNN_train}). This phase yields the optimal set of  parameters learned from the \emph{training dataset}. They are then used to classify new data points from the \emph{testing dataset} during the testing phase of Factional-DNN (\cref{alg:fDNN_test}). We compare the results of our Fractional-DNN with the classical RNN \cref{cont_gen_prob}.

The rest of this section is organized as follows: First, we discuss some data preprocessing and implementation details. Then we describe the datasets being used, and finally we present the experimental results. 

\subsection{Implementation Details}
\begin{enumerate}[label=(\roman*)]
\item \textbf{Batch Normalization.} During the training phase, we use the batch normalization (BN) technique \cite{Szegedy2015_BN}. At each iteration we randomly select a mini-batch, which comprises of $50\%$ of the training data. We then normalize the mini-batch $\hat{Y}_0 \subset Y_0$, to have a zero mean and a standard deviation of one, i.e. 
\begin{equation} \label{eq:BN}
\hat{Y}_0= \frac{\hat{Y}_0 - \mu(\hat{Y}_0)}{s(\hat{Y}_0)},
\end{equation}
where $\mu$ is the mean and $s$ is the standard deviation of the mini-batch. The normalized mini-batch is then used to train the network in that iteration. At the next iteration, a new mini-batch is randomly selected. This process is repeated $m_2$ times. Batch normalization prevents gradient blow-up, helps speed up the learning and reduces the variation in parameters being learned.

Since the design variables are learnt on training data processed with BN, we also process the testing data with BN, in which case the mini-batch is the whole testing data. 

\item \textbf{Activation Function.} For the experiments we have performed, we have used the hyperbolic tangent function as the activation function, for which,
\[
\sigma(x) = \tanh(x), \;\; \text{and} \;\; \sigma^{\prime}(x) = 1-\tanh^2(x).
\]

\item \textbf{Regularization.} In our experiments, we have used the following regularization:
\[
\R(W,K,b) :=\frac{\xi_W}{2}\|W \|^2_F  + \frac{\xi_K}{2N} \|(-\Delta)_hK(t) \|^2_F + \frac{\xi_b}{2N} \| b(t) \|^2_2
\]
where $(-\Delta)_h$ is the discrete Laplacian, and $\xi_W, \xi_K, \xi_b$ are the scalar regularization strengths, and $\| \cdot \|_F$ is the Frobenius norm. 

Notice that with the above regularization, we are enforcing Laplacian-smoothing on $K$. For a more controlled smoothness, one 
could also use the fractional Laplacian regularization introduced in \cite{antil2017spectral}, see also \cite{HAntil_CNRautenberg_2019b} and \cite{ADK_19_ml_tomo}.

\item \textbf{Order of Fractional Time Derivative.} 
In our computations, we have chosen $\gamma$ heuristically. We remark that this fractional exponent on time derivative can be learnt in a similar manner as the fractional exponent on Laplacian was learnt in \cite{ADK_19_ml_tomo}.

\item \textbf{Optimization Solver and Xavier Initialization.} The optimization algorithm we have used is the BFGS method with Armijo line search \cite{Kelly_bfgs}. The stopping tolerance for the BFGS algorithm is set to $1e-6$ or maximum number of optimization iterations $m_2$, whichever is achieved first. However, in our experiments, the latter is achieved first in most cases. The design variables are initialized using Xavier initialization \cite{Bengio2010}, according to which, the biases $b$ are initialized as 0, and the entries of $W$, and $K_j$ are drawn from the uniform distribution $U[-a,a]$. We consider $a=\sqrt{\frac{3}{n_f}}$ for the activation function $\sigma(\cdot) = \tanh(\cdot)$, and $a=\frac{1}{\sqrt{n_f}}$ for other activation functions. 

\item \textbf{Network Layers vs. the Final Time.}
For our experiments, we heuristically choose the number of layers $N$, and the discretization step-length for forward and backward propagation as $\tau=0.2$. Thus our final time is given by, $T=t_N = N\tau$.

\item \textbf{Classification Accuracy.}
We remark that when we calculate $C_{train}=S(W,Y_N)$, we obtain a probability distribution of the samples belonging to the classes. We consider the class with the highest probability as the predicted class. Then, we use a very standard procedure to compare $C_{train}$ with $C_{obs}$.
\[
\text{n}_{\text{cor,train}}:=\text{No. of correctly identified labels} = n - \frac{1}{2} \|C_{obs}-C_{train}\|_F^2.
\]
\[ \text{training error} = 1- \frac{\text{n}_{\text{cor,train}}}{n}, \quad \text{and}  \quad \alpha_{train} = \frac{\text{n}_{\text{cor,train}}}{n} \times 100.\]
The same procedure is used to compute $C_{test}$ and $\alpha_{test}$.

\item \textbf{Gradient Test.}
To verify the gradients in \cref{eq:design_disc}, we perform a gradient test by comparing them to a finite difference gradient approximation of \cref{RNN_frac}. In \cref{f:deriv_test} we show that the two conform and we obtain the expected order of convergence for all the design variables.

\begin{figure}[htb]
\begin{center}
\includegraphics[scale=0.35]{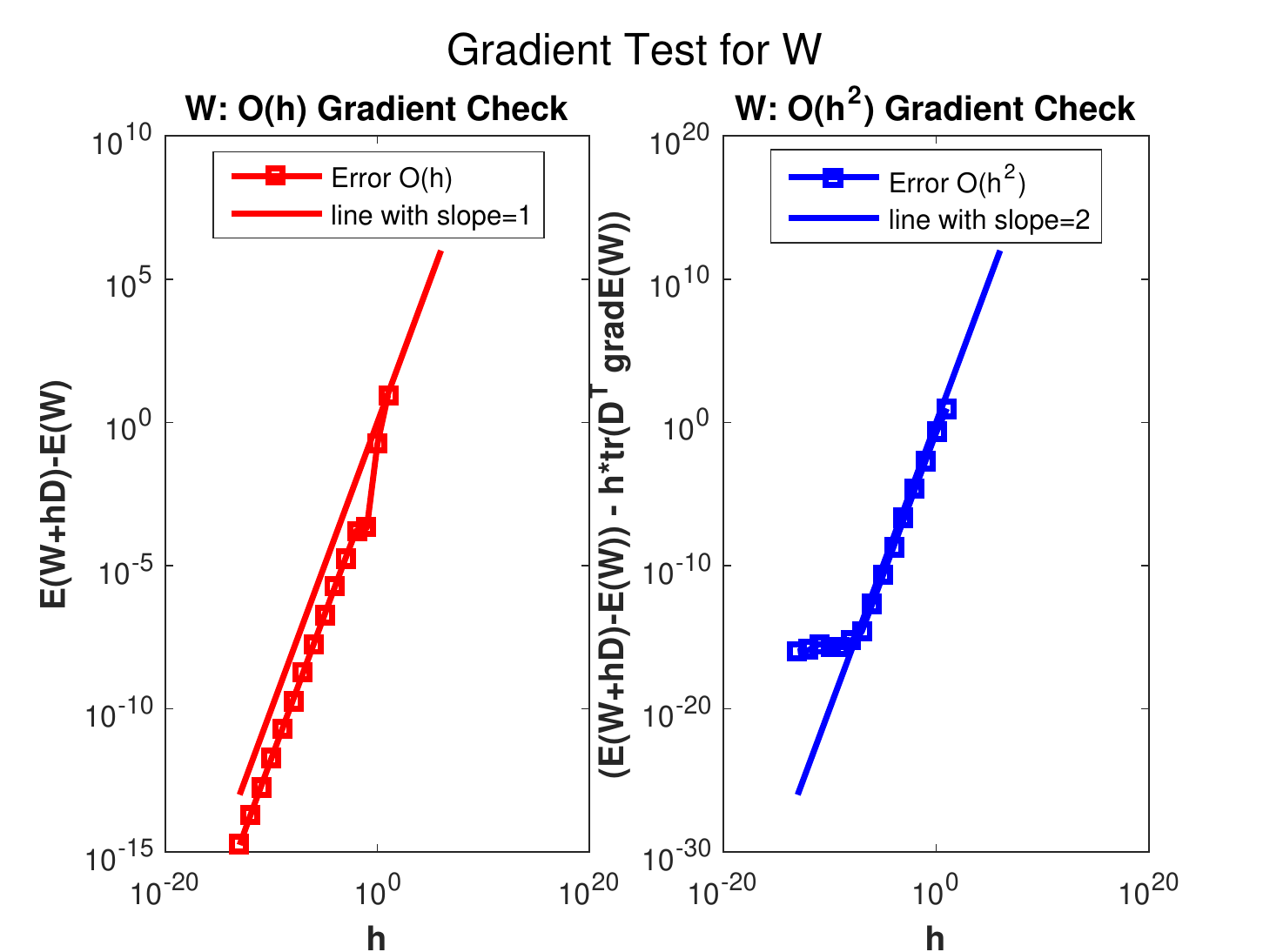}
\includegraphics[scale=0.35]{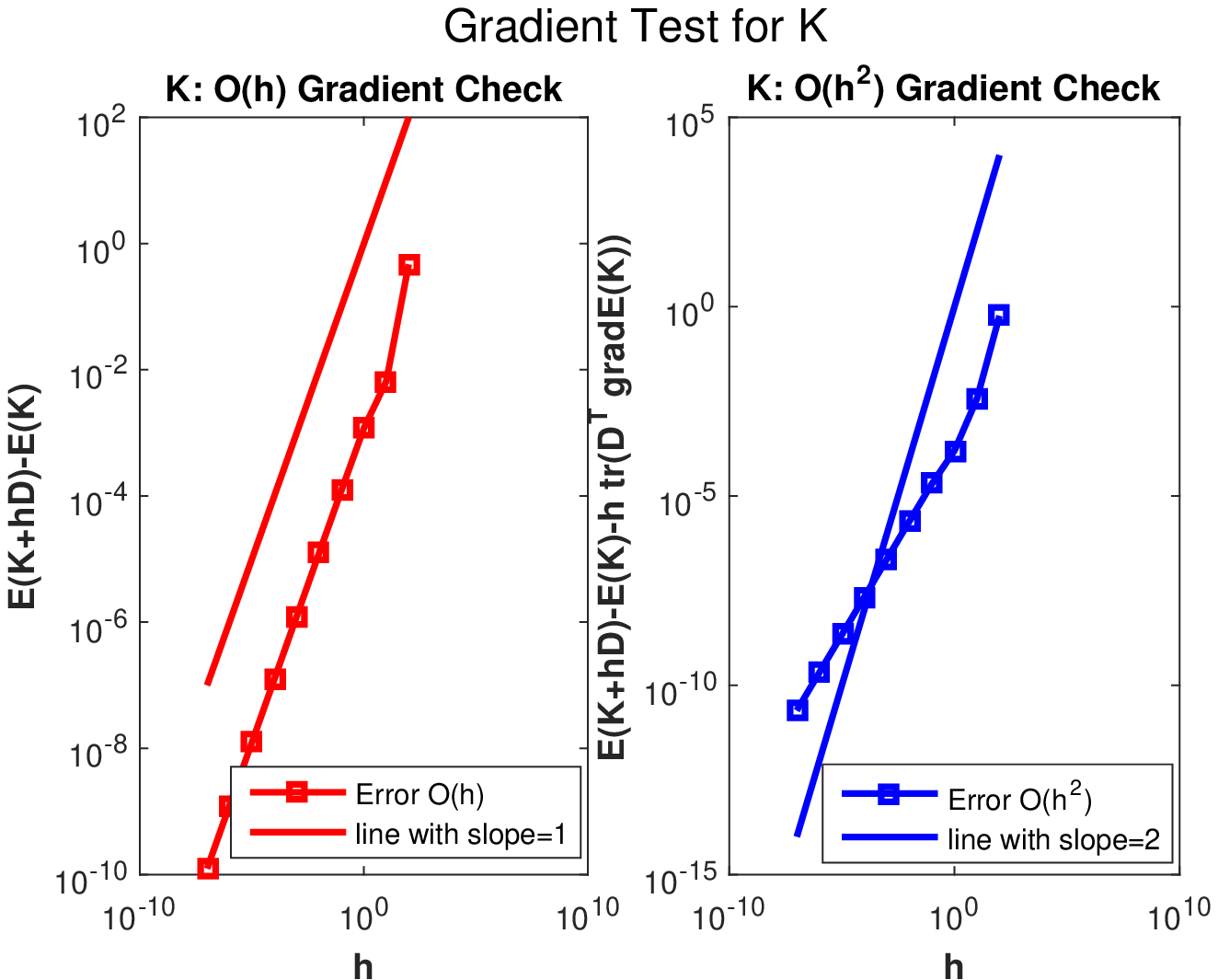}
\includegraphics[scale=0.35]{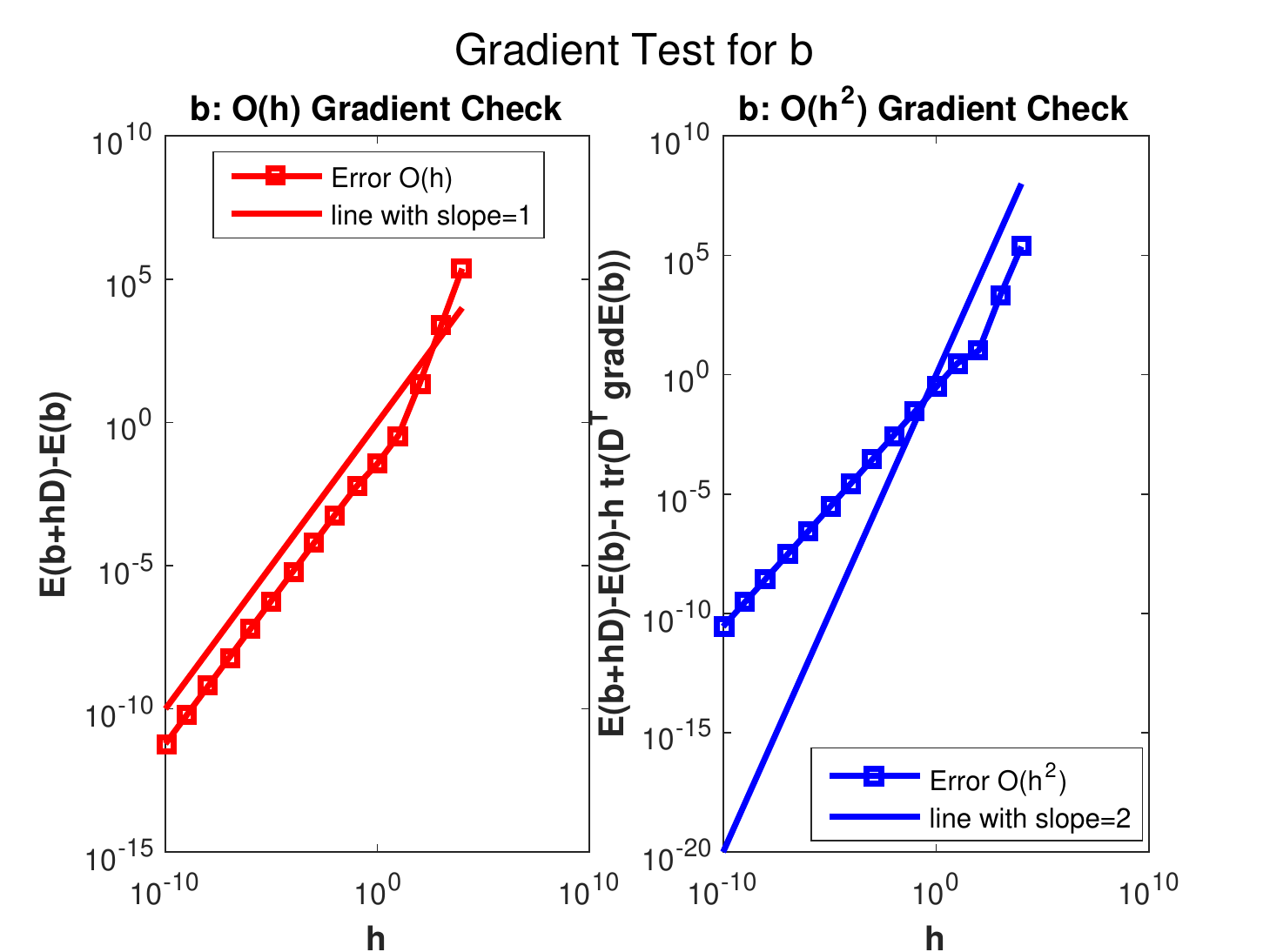}
\caption{\label{f:deriv_test}Comparison between derivative with respect to the design variables 
and finite difference approximation. The expected rate of convergence is obtained.}
\end{center}
\end{figure}

\item \textbf{Computational Platform.} All the computations have been carried out in MATLAB R2015b on a laptop with an Intel Core i7-8550U processor. 
\end{enumerate}
\subsection{Experimental Datasets\label{Datasets}}
We describe the datasets we have used to validate our proposed Fractional-DNN algorithm below.
\begin{enumerate}[label = $\bullet$]
\item \textbf{Dataset 1: Coordinate to Level Set (CLS).} 
This data comprises of a set of 2D coordinates, i.e. $Y_0:= \{ (x_i,y_i)\;\; | \;\; i = 1,\cdots,n;\; (x_i,y_i) \in \mathbb{R}^2([0,1])\}$. Next, we consider the following piecewise function,
\begin{equation} \label{clamped_f}
v(x,y) = \left\{
\begin{aligned}
& 1 \quad &\forall\; x \leq y\\
& 0 &\forall \;x>y
\end{aligned}
\right. \hspace{1cm} \forall \;\; x,y \in [0,1].
\end{equation}
The coordinates are the features in this case, hence $n_f = 2$. Further, we have $n_c=2$ classes, which are the two level sets of $v(x,y)$. Thus, for the $ith$ sample $Y_0^{(i)}, \; C_{obs}^{(i)} \in \RR^{n_c}$ is a standard basis vector which represents the probability that $Y_0^{(i)}$ belongs to the class label $\{1,2\}$.

\item\textbf{Dataset 2: Perfume Data (PD) \cite{Dua:2019, ESME2016452}.}
This dataset comprises of odors of $20$ different perfumes measured via a handheld meter (OMX-GR sensor) every second, for 28 seconds. For this data,  $Y_0:= \{ (x_i,y_i)\;\; | \;\; i = 1,\cdots,n;\; x_i,y_i \in \mathbb{Z}_+\}$, thus $n_f = 2$. The classes, $n_c=20$, pertain to $20$ different perfumes. we construct $C_{obs}$ in the same manner as we did for Dataset 1.
\end{enumerate}

\subsection{Forward Propagation as a Dynamical System}
\begin{figure}[htb]
\begin{center}
\includegraphics[scale=0.55]{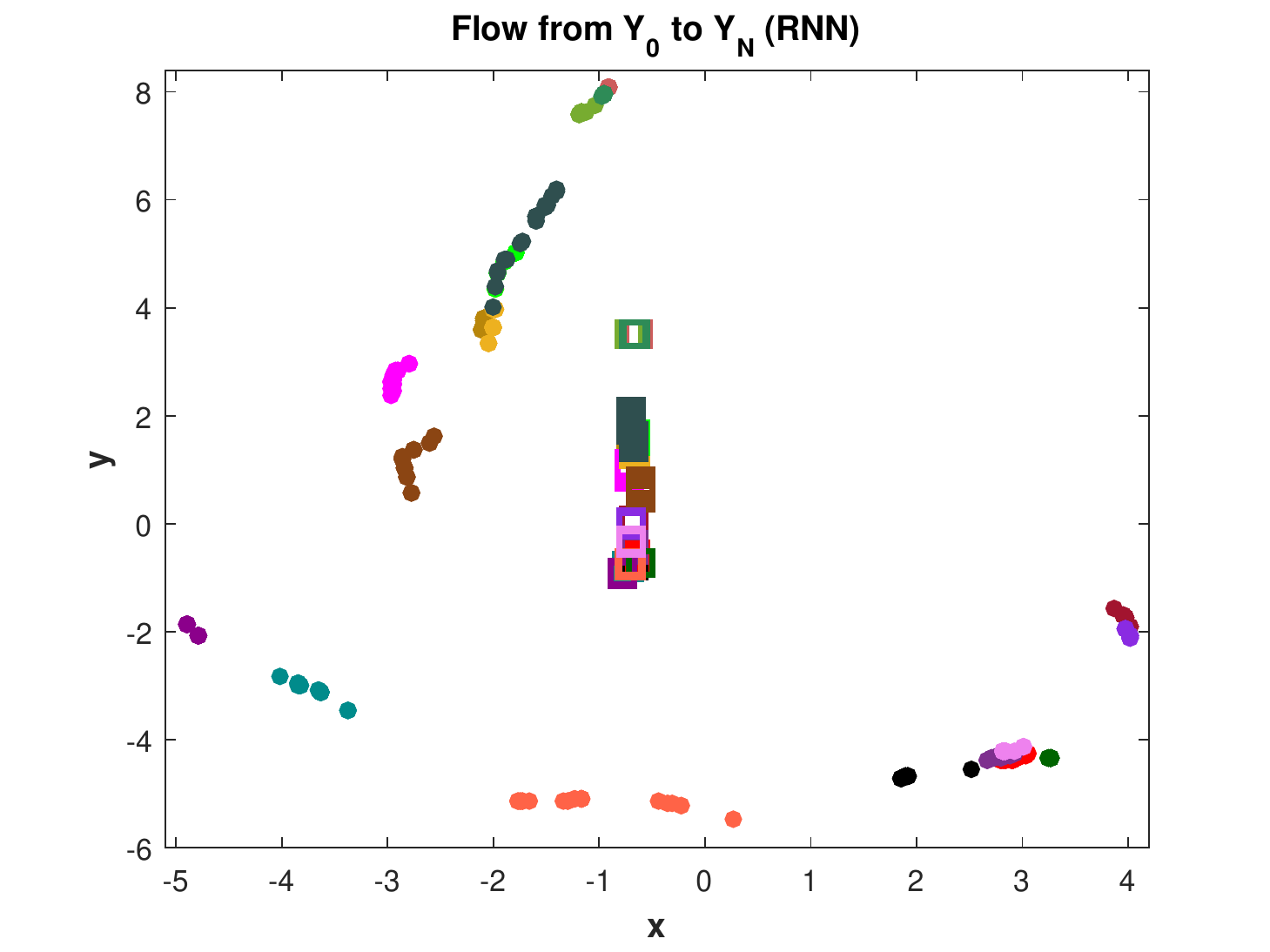}
\includegraphics[scale=0.55]{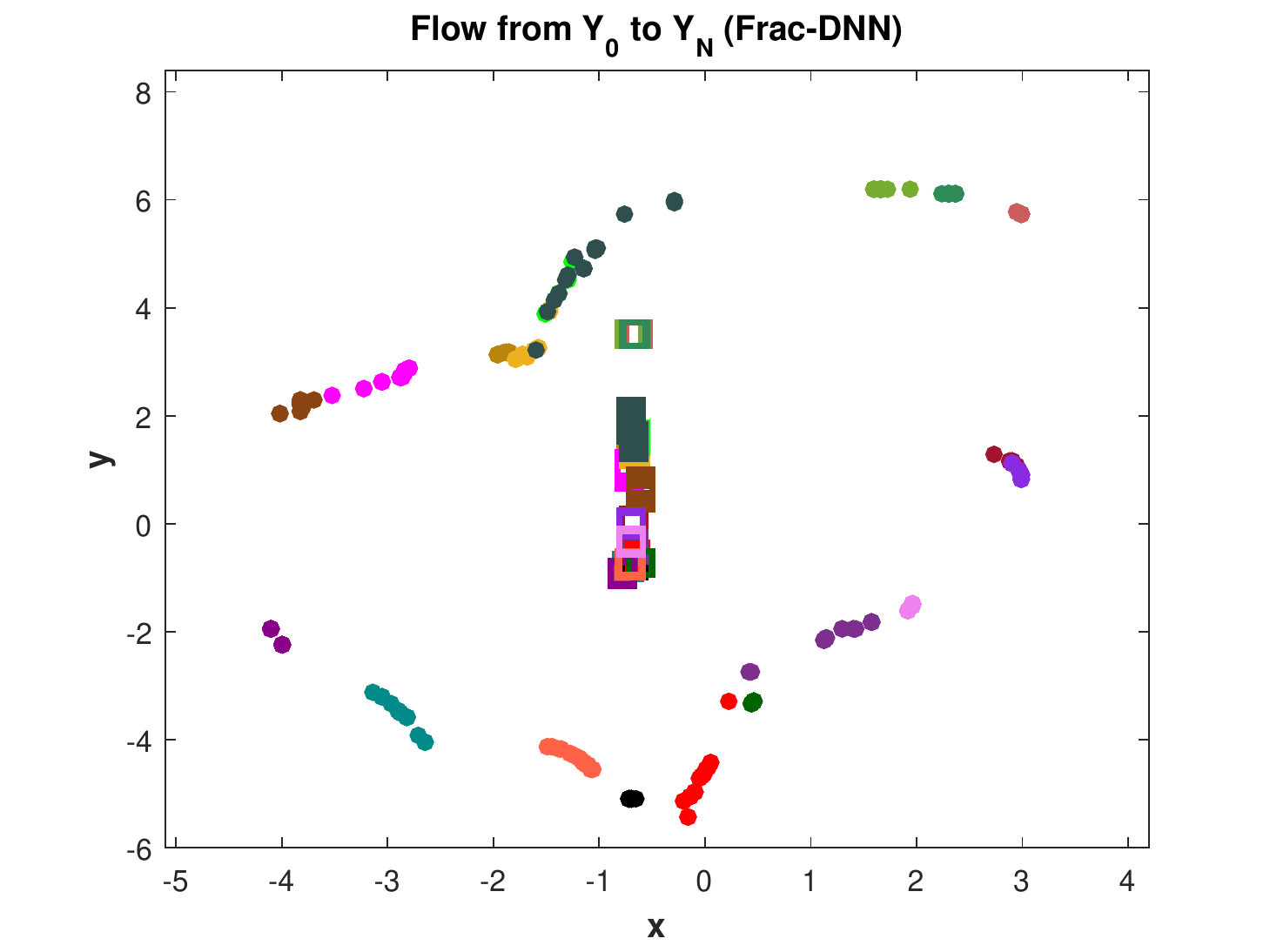}
\caption{\label{fig:DS2} Forward propagation of perfume data from $Y_0$ (\textit{squares}) to $Y_N$ (\textit{dots}) via standard RNN (\textit{left}) and Fractional-DNN (\textit{right}). Note that data is more linearly separable for Fractional-DNN. Different colors represent different classes }
\end{center}
\end{figure}

In the introduction we mentioned the idea of representing a DNN as an optimization problem constrained by a dynamical system. This has turned out to be a strong tool in studying the underlying mathematics of DNNs. In \cref{fig:DS2} we numerically demonstrate how this viewpoint enables a more efficient strategy for distinguishing between the classes. First we consider the perfume data, which has two features, namely the $(x,y)$ coordinates, and let it flow, i.e. forward propagate. When this evolved data is presented to the classifier functional (e.g. softmax function in our case), a spatially well-separated data is easier to classify. We plot the input data $Y_0$, represented as squares, as well as the evolved data $Y_N$ after it has passed through $N$ layers. The $20$ different colors correspond to the $20$ different classes for the data, which help us visually track the evolution from $Y_0$ to $Y_N$. The evolution under standard RNN is shown in the left plot, and that of Fractional-DNN is shown in the right plot. The configuration for these plots is the same as discussed in \cref{sec:exp_results} below and pertains to the trained models. Notice that at the bottom right corner of the RNN evolution plot, the purple, pink and red data points are overlapping which poses a challenge for the classifier to distinguish between the classes. In contrast, Fractional-DNN has separated out those points quite well.

We remark that this separation also gives a hint as to the number of layers needed in a network. We  need enough number of layers which would let the data evolve enough to be easily separable. However, the visualization can get restricted to $n_f \leq3$, therefore for data with $n_f>3$, it may be challenging to get a sense of number of layers needed to make the data separable-enough. 
 
\subsection{Vanishing Gradient Issue}

\begin{figure}[h!] 
\begin{center}
\includegraphics[scale=0.5]{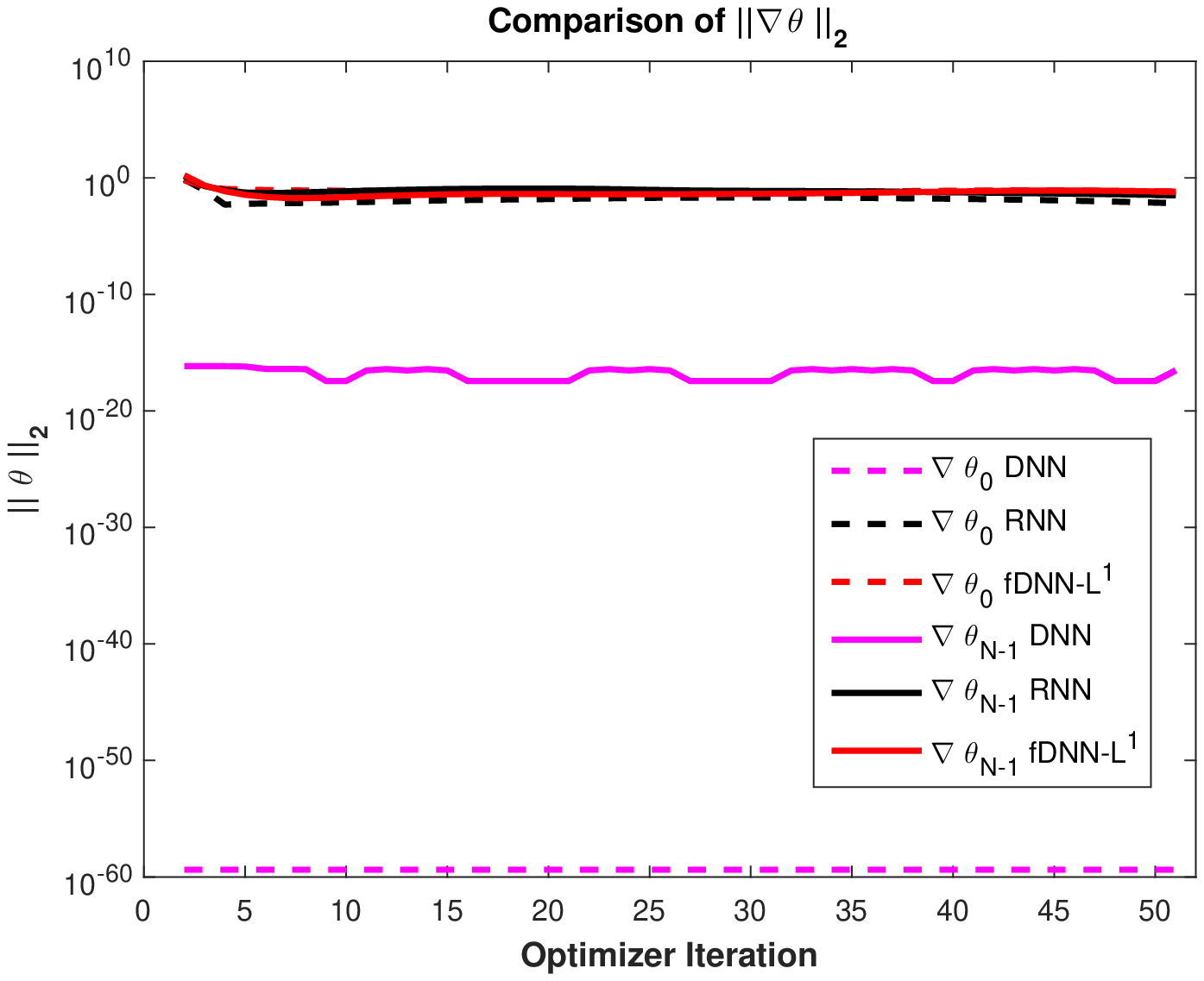} \includegraphics[scale=0.5]{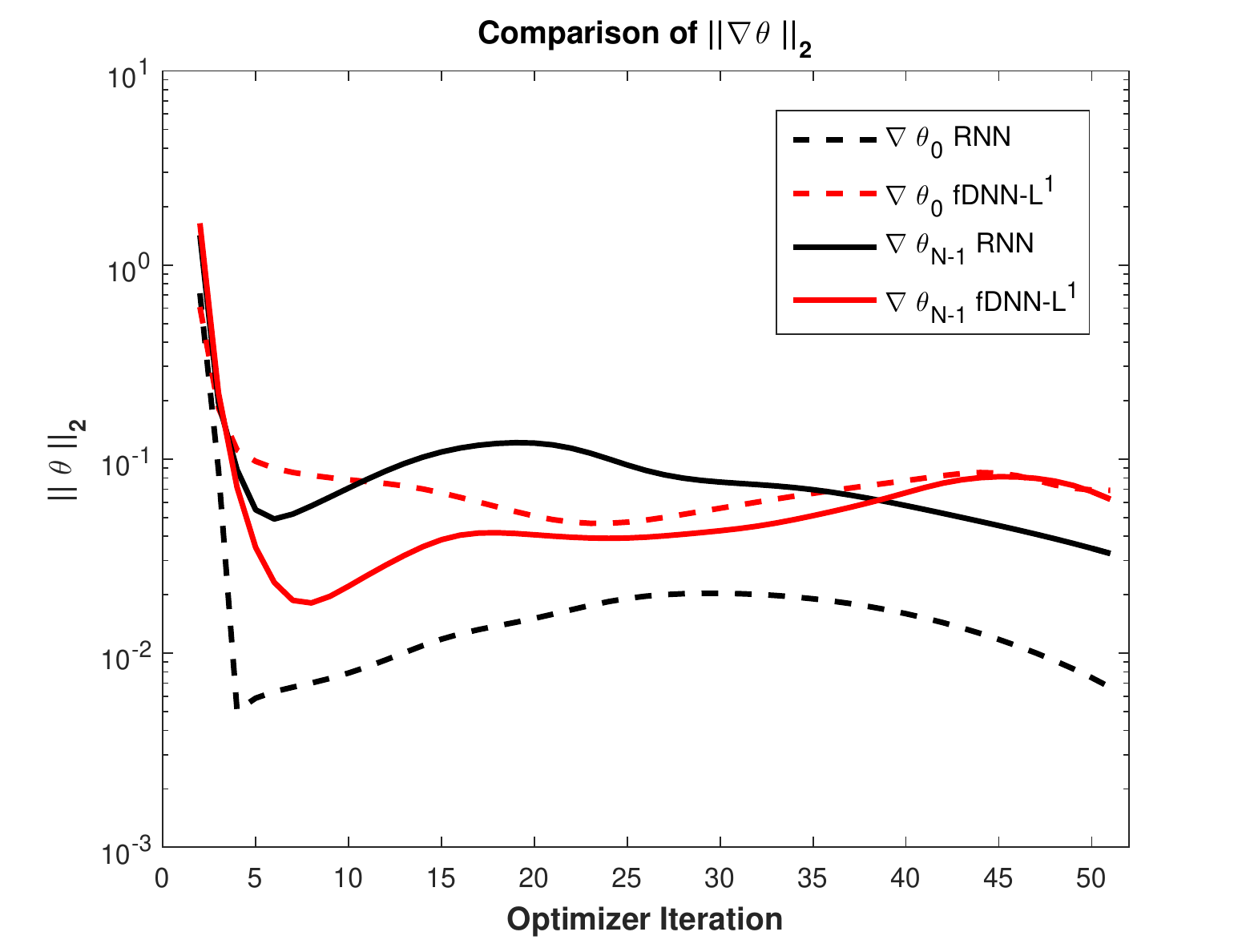} 		\caption{Demonstration of the gradient norm of $\theta = (K,b)$ at the first layer(dotted line) and last layer (solid line) of the network for various algorithms, namely standard DNN (magenta), RNN (black), and Fractional-DNN with $L^1$ scheme approximation (red). The figure on the \textit{right} is the zoomed in version of figure on the \textit{left}. Note the improvement in relative gradient propagation across layers for Fractional-DNN which leads to better learning ability and improves the vanishing gradient issue.}
\label{fig:VG} 
\end{center}
\end{figure}
In earlier sections, we remarked that Fractional-DNN handles the vanishing gradient issue in a better way. The vanishing gradient issue arises when the gradient of the design variables vanishes across the layers as the network undergoes backpropagation, see \cite{VG_DNN_compChem} and references therein. As a consequence, feature extraction in the initial layers gets severely affected, which in turn affects the learning ability of the whole network. We illustrate this phenomenon for the networks under discussion in \Cref{fig:VG}. In the \textit{left} plot of \Cref{fig:VG}, we compare the $\|\cdot\|_2$ of the gradient of design variables $\theta = (K,b)$ against optimization solver (steepest descent in this case) iterations for standard DNN (which does not have any skip connection) in \textit{magenta}, classical RNN  \cref{cont_gen_prob} in \textit{black}, and Fractional-DNN with $L^1$-scheme approximation from \cref{alg:fDNN_train} in red. In the \textit{right} plot of \Cref{fig:VG} we have omitted the standard DNN plot to take a closer look at the other two. Observe that as gradient information propagates backward, i.e. from layer $N-1$ to $0$, its magnitude reduces by one order in the case of standard RNN. This implies that enough information is not being passed to initial layers relative to the last layer. In contrast, the Fractional-DNN is carrying significant information back to the initial layers while maintaining a relative magnitude. This improves the overall health of the network and improves learning. This test has been performed on Perfume Data (Dataset $2$) with $70$ layers and regularization turned off.

\subsection{Experimental Results}\label{sec:exp_results}
We now solve the classification problem \cref{RNN_frac} for the datasets described in \cref{Datasets} via our proposed Fractional-DNN algorithm, presented in \cref{Frac_DNN_alg}. We then compare it with the standard RNN architecture \cref{cont_gen_prob}. The details and results of our experiments are given in \cref{rnn_mem_results}.
\begin{table}[h!] 
\caption{\label{rnn_mem_results}Comparison of classification accuracy for various datasets using the standard RNN \cref{cont_gen_prob} with our proposed Fractional-DNN \cref{cont_gen_prob_frac} with $L^1$ scheme approximation. Note the improvement in results due to Fractional-DNN.}
\begin{center}
\begin{tabular}{|c|c|c|c|c|}
\hline
\textbf{Dataset} 	& \textbf{CLS} 		&\textbf{CLS} 		& \textbf{PD}		& \textbf{PD}\\ \hline
\textbf{Time Derivative}	& Standard   	&Frac-$L^1$			& Standard 	    	&Frac-$L^1$\\ \hline
$\boldsymbol{\gamma}$   & --                 & $0.1$				& --         & $0.9$\\ \hline
\textbf{$\boldsymbol{n_{train}}$}  &$10000$  & $10000$			& $560$ 	   		& $560$\\ \hline
\textbf{$\boldsymbol{n_{test}}$}   & $10000$ & $10000$			& $532$     		& $532$\\ \hline
\textbf{$N$}        &$5$    		& $5$				& $35$              & $35$\\ \hline
\textbf{$\boldsymbol{m_1}$} & $6$       & $6$	& $567$        & $567$\\ \hline
\textbf{$\boldsymbol{m_2}$} & $30$       & $30$	& $15$        & $15$\\ \hline
\textbf{$\boldsymbol{\xi_W}$}          & $1e-1$                    & $1e-1$	& $1e-8$                       & $1e-8$\\ \hline
\textbf{$\boldsymbol{\xi_K}$}          & $1e+2$                          & $1e+2$ 	& $0$             & $0$\\ \hline
\textbf{$\boldsymbol{\xi_b}$}         & $1e-2$                     & $1e-2$ 	& $0$                 & $0$\\ \hline
\textbf{$\boldsymbol{\alpha_{train}}$} & $99.76 \%$                  &$99.82\%$		& $52.86 \%$      &$70.36\%$\\ \hline
\textbf{$\boldsymbol{\alpha_{test}}$}  &$99.79\%$                &$99.79\%$		& $45.49 \%$              &$84.21\%$\\ \hline
\end{tabular}
\end{center}
\end{table}

Note that the results obtained via Fractional-DNN are either comparable to (e.g. for CLS data) or significantly better than (e.g. for PD) the standard RNN architecture.

We remark that while CLS data (Dataset 1) is a relatively simpler problem to solve (two features and two classes), the Perfume Data (Dataset 2) is not. In the latter case, each dataset comprises of only two features, and there are $20$ different classes. Furthermore, the number of available samples for training is small. In this sense, classification of this dataset is a challenging problem. 
There have been some results on classification of perfume data using only the training dataset (divided between training and testing) \cite{ESME2016452}, but to the best of our knowledge, classification on the complete dataset using both the training and testing sets \cite{Dua:2019} is not available.

In our experiments, we have also observed that Fractional-DNN algorithm needs lesser number of Armijo line-search iterations than the standard RNN. This directly reflects an improvement in the learning rate via Fractional-DNN. We remark that in theory, Fractional-DNN should use memory more efficiently than other networks, as it encourages feature reuse in the network. 

\section{Discussion\label{Disc}}

There is a growing body of research which indicates that deep learning algorithms, e.g. a residual neural network, can be cast as optimization problems constrained by ODEs or PDEs. In addition, thinking of continuous optimization problems can make the approaches machine/architecture independent. This opens a plethora of tools from constrained optimization theory which can be used to study, analyze, and enhance the deep learning algorithms. Currently, the mathematical foundations of many machine learning models are largely lacking. Their success is mostly attributed to empirical evidence. Hence, due to the lack of mathematical foundation, it becomes challenging to fix issues, like network instability, vanishing and exploding gradients, long training times, inability to approximate non-smooth functions, etc., when a network breaks down. 

In this work we have developed a novel continuous model and stable discretization of deep neural networks that incorporate history. In particular, we have developed a fractional deep neural network (Fractional-DNN) which allows the network to admit memory across all the subsequent layers. We have established this via an optimal control problem formulation of a deep  neural network bestowed with a fractional time Caputo derivative. We have then derived the optimality conditions using the Lagrangian formulation. We have also discussed discretization of the fractional time Caputo derivative using $L^1$-scheme and presented the algorithmic framework for the discretization. 

We expect that by keeping track of history in this manner improves the vanishing gradient problem and can potentially  strengthen feature propagation, encourage feature reuse and reduce the number of unknown parameters. We have numerically illustrated the improvement in the vanishing gradient issue via our proposed Fractional-DNN. We have shown that Fractional-DNN is better capable of passing information across the network layers which maintains the relative gradient magnitude across the layers, compared to the standard DNN and standard RNN. This allows for a more meaningful feature extraction to happen at each layer. 

We have shown successful application of Fractional-DNN for classification problems using various datasets, namely the Coordinate to Level Set (CLS dataset) and Perfume Data. We have compared the results against the standard-RNN and have shown that the Fractional-DNN algorithm yields improved results. 

We emphasis that our proposed Fractional-DNN architecture has a \emph{memory} effect due to the fact that it allows propagation of features in a cumulative manner, i.e. at each layer all the precedent layers are visible. Reusing the network features in this manner reduces the number of parameters that the network needs to learns in each subsequent layer. Fractional-DNN has a rigorous mathematical foundation and algorithmic framework which establishes a deeper understanding of deep neural networks with memory. This enhances their applicability to scientific and engineering applications.

We remark that code optimization is part of our forthcoming work. This would involve efficient Graphic Processing Unit usage and parallel computing capabilities. We also intend to develop a python version of this code and incorporate it into popular deep learning libraries like TensorFlow, PyTorch etc. We are also interested in expanding the efficiency of this algorithm to large-scale problems suitable for High Performance Computing.

\section{Acknowledgments}

The authors would like to thank Prasanna Balaprakash, Tamara G. Kolda and Lars Ruthotto for several discussions and comments during the course of this project. 

\bibliographystyle{plain}
\bibliography{ref_fracDNN}

\begin{thebibliography}{10}

\bibitem{anderson1965iterative}
D.~G. Anderson.
\newblock Iterative procedures for nonlinear integral equations.
\newblock {\em Journal of the ACM (JACM)}, 12(4):547--560, 1965.

\bibitem{antil2017spectral}
H.~Antil and S.~Bartels.
\newblock Spectral {A}pproximation of {F}ractional {PDE}s in {I}mage
  {P}rocessing and {P}hase {F}ield {M}odeling.
\newblock {\em Comput. Methods Appl. Math.}, 17(4):661--678, 2017.

\bibitem{ADK_19_ml_tomo}
H.~Antil, Z.~Di, and R.~Khatri.
\newblock Bilevel optimization, deep learning and fractional laplacian
  regularization with applications in tomography.
\newblock {\em Inverse Problems}, 2020.

\bibitem{antil_frac_time}
H.~Antil, C.~Lizama, R.~Ponce, and M.~Warma.
\newblock Convergence of solutions of discrete semi-linear space-time
  fractional evolution equations.
\newblock {\em arXiv preprint arXiv:1910.07358}, 2019.

\bibitem{antil_otarola_salgado}
H.~Antil, E.~Ot\'{a}rola, and A.J. Salgado.
\newblock A space-time fractional optimal control problem: analysis and
  discretization.
\newblock {\em SIAM J. Control Optim.}, 54(3):1295--1328, 2016.

\bibitem{HAntil_CNRautenberg_2019b}
H.~Antil and C.N. Rautenberg.
\newblock Sobolev spaces with non-{M}uckenhoupt weights, fractional elliptic
  operators, and applications.
\newblock {\em SIAM J. Math. Anal.}, 51(3):2479--2503, 2019.

\bibitem{Bengio2012LearningRate}
Y.~Bengio.
\newblock {\em Practical Recommendations for Gradient-Based Training of Deep
  Architectures}, pages 437--478.
\newblock Springer Berlin Heidelberg, Berlin, Heidelberg, 2012.

\bibitem{Bengio1994}
Y.~Bengio, P.~Simard, and P.~Frasconi.
\newblock Learning long-term dependencies with gradient descent is difficult.
\newblock {\em IEEE Transactions on Neural Networks}, 5(2):157,166, 1994-03.

\bibitem{benning2019deep}
M.~Benning, E.~Celledoni, M.~Ehrhardt, B.~Owren, and C.-B. Schönlieb.
\newblock Deep learning as optimal control problems: Models and numerical
  methods.
\newblock {\em Journal of Computational Dynamics}, 6:171--198, 01 2019.

\bibitem{bischke2017detection}
B.~Bischke, P.~Bhardwaj, A.~Gautam, P.~Helber, D.~Borth, and A.~Dengel.
\newblock Detection of flooding events in social multimedia and satellite
  imagery using deep neural networks.
\newblock In {\em Working Notes Proceedings of the MediaEval 2017. MediaEval
  Benchmark, September 13-15, Dublin, Ireland}. MediaEval, 2017.

\bibitem{brown2018analysis}
T.~Brown, S.~Du, H.~Eruslu, and F.-J. Sayas.
\newblock Analysis of models for viscoelastic wave propagation.
\newblock {\em arXiv preprint arXiv:1802.00825}, 2018.

\bibitem{Haber2017_RNN}
B.~Chang, L.~Meng, E.~Haber, F.~Tung, and D.~Begert.
\newblock Multi-level residual networks from dynamical systems view.
\newblock {\em arXiv preprint arXiv:1710.10348}, 2017.

\bibitem{chen2018voxresnet_bio}
H.~Chen, Q.~Dou, L.~Yu, J.~Qin, and P.-A. Heng.
\newblock Voxresnet: Deep voxelwise residual networks for brain segmentation
  from 3d mr images.
\newblock {\em NeuroImage}, 170:446--455, 2018.

\bibitem{ResNetPlus2018}
K.~{Chen}, K.~{Chen}, Q.~{Wang}, Z.~{He}, J.~{Hu}, and J.~{He}.
\newblock Short-term load forecasting with deep residual networks.
\newblock {\em IEEE Transactions on Smart Grid}, 10(4):3943--3952, July 2019.

\bibitem{CIRESAN2012333}
D.~Cire\c{s}an, U.~Meier, J.~Masci, and J.~Schmidhuber.
\newblock Multi-column deep neural network for traffic sign classification.
\newblock {\em Neural Networks}, 32:333 -- 338, 2012.
\newblock Selected Papers from IJCNN 2011.

\bibitem{cortes2016adanet}
C.~Cortes, X.~Gonzalvo, V.~Kuznetsov, M.~Mohri, and S.~Yang.
\newblock Adanet: Adaptive structural learning of artificial neural networks.
\newblock {\em Efficient Methods for Deep Neural Networks (EMDNN)}, 07 2016.

\bibitem{Dua:2019}
D.~Dua and C.~Graff.
\newblock {UCI} machine learning repository, 2017.

\bibitem{Weinan_2019}
W.~E.
\newblock Machine learning: Mathematical theory and scientific applications.
\newblock {\em Notices of the American Mathematical Society},
  66(11):1813–1820, Dec 2019.

\bibitem{ESME2016452}
E.~Esme and B.~Karlik.
\newblock Fuzzy c-means based support vector machines classifier for perfume
  recognition.
\newblock {\em Applied Soft Computing}, 46:452 -- 458, 2016.

\bibitem{Bengio2010}
X.~Glorot and Y.~Bengio.
\newblock Understanding the difficulty of training deep feedforward neural
  networks.
\newblock In {\em Proceedings of the Thirteenth International Conference on
  Artificial Intelligence and Statistics}, volume~9 of {\em Proceedings of
  Machine Learning Research}, pages 249--256, Chia Laguna Resort, Sardinia,
  Italy, 13--15 May 2010. PMLR.

\bibitem{VG_DNN_compChem}
G.~B. Goh, N.~O. Hodas, and A.~Vishnu.
\newblock Deep learning for computational chemistry.
\newblock {\em Journal of Computational Chemistry}, 38(16):1291--1307, 2017.

\bibitem{goldt2019modelling}
S.~Goldt, M.~M{\'e}zard, F.~Krzakala, and L.~Zdeborov{\'a}.
\newblock Modelling the influence of data structure on learning in neural
  networks.
\newblock {\em arXiv preprint arXiv:1909.11500}, 2019.

\bibitem{gunther2018layer}
S.~Günther, L.~Ruthotto, J.~Schroder, E.~Cyr, and N.~Gauger.
\newblock Layer-parallel training of deep residual neural networks.
\newblock {\em SIAM Journal on Mathematics of Data Science}, 2:1--23, 01 2020.

\bibitem{haber2017stable}
E.~Haber and L.~Ruthotto.
\newblock Stable architectures for deep neural networks.
\newblock {\em Inverse Problems}, 34(1):014004, 2017.

\bibitem{Pock2017_MRI}
K.~Hammernik, T.~Klatzer, E.~Kobler, M.~P. Recht, D.~K. Sodickson, T.~Pock, and
  F.~Knoll.
\newblock Learning a variational network for reconstruction of accelerated mri
  data.
\newblock {\em Magnetic Resonance in Medicine}, 79(6):3055--3071, 6 2018.

\bibitem{he2016deep}
K.~He, X.~Zhang, S.~Ren, and J.~Sun.
\newblock Deep residual learning for image recognition.
\newblock In {\em Proceedings of the IEEE Conference on Computer Vision and
  Pattern Recognition}, pages 770--778, 2016.

\bibitem{huang2017denseNet}
G.~Huang, Z.~Liu, L.~van~der Maaten, and K.~Q. Weinberger.
\newblock Densely connected convolutional networks.
\newblock In {\em Proceedings of the IEEE Conference on Computer Vision and
  Pattern Recognition}, pages 4700--4708, 2017.

\bibitem{Imaizumi2018_nonsmooth}
M.~Imaizumi and K.~Fukumizu.
\newblock Deep neural networks learn non-smooth functions effectively.
\newblock {\em arXiv preprint arXiv:1802.04474}, 2018.

\bibitem{Szegedy2015_BN}
S.~Ioffe and C.~Szegedy.
\newblock Batch normalization: Accelerating deep network training by reducing
  internal covariate shift.
\newblock In {\em Proceedings of the 32nd International Conference on
  International Conference on Machine Learning - Volume 37}, ICML’15, page
  448–456. JMLR.org, 2015.

\bibitem{CNN_tomo_2017}
K.~H. {Jin}, M.~T. {McCann}, E.~{Froustey}, and M.~{Unser}.
\newblock Deep convolutional neural network for inverse problems in imaging.
\newblock {\em IEEE Transactions on Image Processing}, 26(9):4509--4522, Sep.
  2017.

\bibitem{Kelly_bfgs}
C.~T. Kelley.
\newblock {\em Iterative methods for optimization}.
\newblock Frontiers in applied mathematics. SIAM, Philadelphia, 1999.

\bibitem{kilbas_srivastava_fde}
A.~A. Kilbas, H.~M. Srivastava, and J.~J. Trujillo.
\newblock {\em Theory and applications of fractional differential equations},
  volume 204 of {\em North-Holland Mathematics Studies}.
\newblock Elsevier Science B.V., Amsterdam, 2006.

\bibitem{lee2018deepres_bio}
D.~Lee, J.~Yoo, S.~Tak, and J.~C. Ye.
\newblock Deep residual learning for accelerated mri using magnitude and phase
  networks.
\newblock {\em IEEE Transactions on Biomedical Engineering}, 65(9):1985--1995,
  2018.

\bibitem{lu2017finite}
Y.~Lu, A.~Zhong, Q.~Li, and B.~Dong.
\newblock Beyond finite layer neural networks: Bridging deep architectures and
  numerical differential equations, 2017.

\bibitem{mallat2013deep}
S.~Mallat and I.~Waldspurger.
\newblock Deep learning by scattering.
\newblock {\em arXiv preprint arXiv:1306.5532}, 2013.

\bibitem{metzler2000random}
R.~Metzler and J.~Klafter.
\newblock The random walk's guide to anomalous diffusion: a fractional dynamics
  approach.
\newblock {\em Physics reports}, 339(1):1--77, 2000.

\bibitem{NocedalWright2006}
J.~Nocedal and S.~Wright.
\newblock {\em {Numerical Optimization}}.
\newblock Springer Series in Operations Research and Financial Engineering.
  Springer Science {\&} Business Media, New York, December 2006.

\bibitem{Podlubny1999}
I.~Podlubny.
\newblock {\em Fractional differential equations: an introduction to fractional
  derivatives, fractional differential equations, to methods of their solution
  and some of their applications}.
\newblock Mathematics in Science and Engineering. Academic Press, London, 1999.

\bibitem{Qiu2016}
J.~Qiu, Q.~Wu, G.~Ding, Y.~Xu, and S.~Feng.
\newblock A survey of machine learning for big data processing.
\newblock {\em EURASIP Journal on Advances in Signal Processing}, 2016(1):67,
  May 2016.

\bibitem{roux2012stochastic}
N.~L. Roux, M.~Schmidt, and F.~R. Bach.
\newblock A stochastic gradient method with an exponential convergence rate for
  finite training sets.
\newblock In {\em Advances in Neural Information Processing Systems}, pages
  2663--2671, 2012.

\bibitem{ruthotto2018deep}
L.~Ruthotto and E.~Haber.
\newblock Deep neural networks motivated by partial differential equations.
\newblock {\em Journal of Mathematical Imaging and Vision}, 2019.

\bibitem{Samko1993}
S.~G. Samko, A.~A. Kilbas, and O.~I. Marichev.
\newblock {\em Fractional integrals and derivatives}.
\newblock Gordon and Breach Science Publishers, Yverdon, 1993.

\bibitem{scherer2020}
R.~Scherer.
\newblock {\em Computer vision methods for fast image classification and
  retrieval}.
\newblock Springer, 2020.

\bibitem{Carola_resnet_2019}
CB. Sch\"oenlieb, M.~Benning, M.~Ehrhardt, B.~Owren, and E.~Celledoni.
\newblock Research data supporting ``deep learning as optimal control
  problems".
\newblock Dataset, 2019.

\bibitem{Srivastava2015HighwayNet}
R.~K. Srivastava, K.~Greff, and J.~Schmidhuber.
\newblock Training very deep networks.
\newblock In C.~Cortes, N.~D. Lawrence, D.~D. Lee, M.~Sugiyama, and R.~Garnett,
  editors, {\em Advances in Neural Information Processing Systems 28}, pages
  2377--2385. Curran Associates, Inc., 2015.

\bibitem{Mstynes}
M.~Stynes.
\newblock Too much regularity may force too much uniqueness.
\newblock {\em Fract. Calc. Appl. Anal.}, 19(6):1554--1562, 2016.

\bibitem{Tai_2017_CVPR}
Y.~{Tai}, J.~{Yang}, and X.~{Liu}.
\newblock Image super-resolution via deep recursive residual network.
\newblock In {\em 2017 IEEE Conference on Computer Vision and Pattern
  Recognition (CVPR)}, pages 2790--2798, July 2017.

\bibitem{veit2016residual}
A.~Veit, M.~J. Wilber, and S.~Belongie.
\newblock Residual networks behave like ensembles of relatively shallow
  networks.
\newblock In {\em Advances in neural information processing systems}, pages
  550--558, 2016.

\bibitem{Wigderson2019}
A.~Wigderson.
\newblock {\em Mathematics and Computation}.
\newblock A Theory Revolutionizing Technology and Science. Princeton University
  Press, Princeton and Oxford, 2019.
\newblock (in press).

\bibitem{wu2018deep_imaging}
S.~Wu, S.~Z., and Y.~L.
\newblock Deep residual learning for image steganalysis.
\newblock {\em Multimedia tools and applications}, 77(9):10437--10453, 2018.

\bibitem{zhang2018missing}
Q.~Zhang, Q.~Yuan, C.~Zeng, X.~Li, and Y.~Wei.
\newblock Missing data reconstruction in remote sensing image with a unified
  spatial--temporal--spectral deep convolutional neural network.
\newblock {\em IEEE Transactions on Geoscience and Remote Sensing},
  56(8):4274--4288, 2018.

\bibitem{Zhang_2018}
Y.~Zhang, Y.~Tian, Y.~Kong, B.~Zhong, and Y.~Fu.
\newblock Residual dense network for image super-resolution.
\newblock {\em 2018 IEEE/CVF Conference on Computer Vision and Pattern
  Recognition}, June 2018.

\end{thebibliography}
\end{document}